\def\sphere{\overline {\Bbb C}}
\def\dist{\text{dist}}
\def\diam{\text{diam}}
\def\dim{\text{dim}}

\documentstyle[12pt]{amsart}
\theoremstyle{plain}                    
\newtheorem{thm}{Theorem}[section]
\newtheorem{cor}[thm]{Corollary}
\newtheorem{lemma}[thm]{Lemma}
\newtheorem{prop}[thm]{Proposition}

\theoremstyle{definition}               

\theoremstyle{remark}                   

\numberwithin{equation}{section}
\numberwithin{figure}{section}

    {

     \goodbreak\smallskip
     \centerline{\tt#1}
     \smallskip
     \nopagebreak[4]
     \begin{itemize}
            \sf\small%
            \baselineskip=1.1\baselineskip %
            }%
    {\end{itemize}\medskip}

\newcommand{\thmref}[1]{Theorem~\ref{#1}}
\newcommand{\propref}[1]{Proposition~\ref{#1}}
\newcommand{\corref}[1]{Corollary~\ref{#1}}
\newcommand{\lemref}[1]{Lemma~\ref{#1}}

\newcommand{\secref}[1]{Section~\ref{#1}}

\newcommand{\C}{{\Bbb C}}
\newcommand{\Cbar}{\overline{\C}}

\advance\marginparwidth by -\textwidth          
\advance\marginparwidth by -1 true in           
\advance\marginparwidth by -\oddsidemargin      
\advance\marginparwidth by -\marginparsep       
\advance\marginparwidth by -.25 true in         
\begin{document}
\baselineskip=24pt
\renewcommand{\thepage}{}              
\title  [Hausdorff dimension and Kleinian groups]
        {Hausdorff dimension and Kleinian groups}

\subjclass{30F40 ; Secondary 28A78, 30C35, 31A15  }

\keywords{ Harmonic measure, Hausdorff dimension, Kleinian 
groups, limit sets, degenerate groups, geometrically infinite,
Teichm{\"u}ller space,
 conical limit set, Poincar{\'e} exponent}

\author {Christopher J. Bishop}
\address{C.J. Bishop\\
         Mathematics Department\\
         SUNYSB \\
         Stony Brook, NY 11794}
\email {bishop@@math.sunysb.edu}
\author {Peter W. Jones}
\address{P.W. Jones\\
         Mathematics Department\\
         Yale University\\
         New Haven, CT 06520 }
\thanks{The first author is partially supported by NSF Grant DMS 92-04092 and
        an Alfred P. Sloan research fellowship. The second author is 
        partially supported by NSF grant DMS-92-13595}

\date{April 4, 1994}
\maketitle

\begin{abstract}
Let $G$ be a non-elementary, finitely generated Kleinian group, $\Lambda(G)$
its limit set and $\Omega(G) = \sphere \backslash \Lambda(G)$
its set of discontinuity. Let $\delta(G)$ be the critical 
exponent for the Poincar{\'e} series and let $\Lambda_c$ be 
the conical limit set of $G$. 
We prove that 
\begin{enumerate}
\item $\delta(G) = \dim(\Lambda_c)$.
\item A simply connected component  $\Omega$
is either a disk or $\dim(\partial \Omega)>1$.
\item $ \Lambda(G)$ is  either totally 
disconnected,  a circle or has
dimension $>1$,
\item  $G$ is geometrically infinite iff $\dim(\Lambda)=2$.
\item If $G_n \to G$ algebraically then $\dim(\Lambda)\leq \liminf
\dim(\Lambda_n)$.
\item The Minkowski dimension of $\Lambda$  equals the Hausdorff 
dimension.
\item If $\text{area}(\Lambda)=0$ then $\delta(G) =\dim(\Lambda(G))$.
\end{enumerate}

\end{abstract}

%
\clearpage
\setcounter{page}{1}
\renewcommand{\thepage}{\arabic{page}}
\section{Statement of results}  \label{statement}

Consider a group $G$ of M{\"o}bius transformations acting 
on the two sphere $S^2$.
Such a group $G$ also acts on the hyperbolic $3$-ball 
$\Bbb B$ (with its boundary $S^2$ identified with 
the Riemann sphere).
The limit set, $\Lambda(G)$, is the accumulation set 
(on $S^2$) of the orbit of the origin. We say the 
group is discrete if is discrete as a subgroup of 
$\text{PSL}(2, \Bbb C)$, (i.e., if the identity element 
is isolated). 
The ordinary set of  $G$, $\Omega(G)$, is the subset of 
$S^2$ where $G$ acts discontinuously, e.g., $\Omega(G)$ is the 
set of points $z$ such that there exists a disk around $z$ which 
hits itself only finitely often under the action of $G$. If 
$G$ is discrete, then $\Omega(G)= S^2 \setminus \Lambda(G)$.
$G$ is called a Kleinian if it is discrete and $\Omega(G)$ is 
non-empty.  The limit set $\Lambda(G)$ has either 0,1,2 or 
infinitely many points and $G$ is called  elementary if 
$\Lambda(G)$ is finite.

 If we form the Poincar{\'e} series for a discrete $G$ on ${\Bbb B}^3$, 
$$ \sum_G \exp (-s \rho(0,g(0))),$$
($\rho$ is the hyperbolic metric in ${\Bbb B}^3$) the series
converges for all $s > s_0$ and the minimal such $s_0$ is 
called the {\it Poincar{\'e} exponent} of the group and is denoted
$\delta(G)$.
 A point $ x \in \Lambda(G)$ is called a
{\it conical limit} point if there is a sequence of orbit points
which converges to $x$ inside a (Euclidean) non-tangential cone with 
vertex at $x$ (such points are also called radial limit
points or points of approximation). The set of such points
is denoted $\Lambda_c(G)$. $G$ is called {\it geometrically 
finite} if there is a finite sided fundamental polyhedron 
for the $G$ action on $\Bbb B$. However, the only property 
we shall use is a result of Beardon and Maskit \cite{Beardon-Maskit74}
that $G$ is 
geometrically finite iff $\Lambda(G)$ is the union of 
$\Lambda_c(G)$, the rank 2 parabolic cusps and doubly cusped rank
1 parabolic fixed points of $G$.  This makes it clear that
$\dim(\Lambda_c) = \dim(\Lambda)$ and $\text{area}(\Lambda)=0$
in the geometrically finite case.
(For our purposes we could take the Beardon-Maskit characterization 
to be the definition of geometrically finite.)
By Selberg's lemma any finitely generated discrete group contains
a finite index subgroup without torsion. This subgroup must have
the same limit set as the original group, so for our purposes it 
will always be sufficient to assume $G$ has no torsion.

In this paper ``circle''  will always refer to 
Euclidean circles or lines (e.g., circles on the sphere)
as opposed to topological circles which we will always call 
closed Jordan curves.  Similarly  the terms ``disk'' or ``ball''
will always denote Euclidean balls.

The principle objective of this 
paper is to establish the following 
results.

\begin{thm}\label{delta=conical1}
If $G$ is a  non-elementary, discrete M{\"o}bius  group on 
$\Bbb B$  then $\delta(G) = \dim(\Lambda_c(G))$.
\end{thm}

\begin{thm} \label{thm1}
Suppose $G$ is a finitely generated Kleinian group
and  $\Omega$ is a simply connected, invariant component of $\Omega(G)$.
Then $\dim(\partial \Omega)=1$ iff   
 $\partial \Omega$ is a  circle.
\end{thm}

\begin{thm} \label{thm2}
If $G$ is a finitely generated Kleinian group then its 
limit set is either totally disconnected, a circle or has 
Hausdorff dimension $>1$.
\end{thm}

\begin{thm} \label{gi->dim=2}
If $G$ is a finitely generated Kleinian group  which is geometrically infinite 
then $\dim(\Lambda(G)) =2$.
\end{thm}

\begin{thm}\label{Hdim=lsc}
If $\{G_n\}$ is a sequence of $N$-generated Kleinian groups which converges
algebraically to $G$ then $\dim(\Lambda(G)) \leq \liminf_n 
\dim(\Lambda(G_n))$.
\end{thm}

\begin{thm} \label{H=M1} 
If $G$ is a finitely generated Kleinian group then the
Minkowski dimension of $\Lambda$ exists and equals the
Hausdorff dimension.
\end{thm}

\begin{thm} \label{delta=dim1}
Suppose that  $G$ is a non-elementary, finitely generated Kleinian group and 
that $\text{area}(\Lambda(G)) =0$. Then $\delta(G) = \dim(\Lambda)$.
\end{thm}

These results are arranged roughly in order of dependence and
the amount of required prerequisites from the theory of Kleinian 
groups and hyperbolic manifolds. 
\thmref{delta=conical1} uses nothing
but the definitions and a few simple properties of M{\"o}bius 
transformations.
We have only been able to locate this result in the 
literature under the additional assumptions that 
$G$ is geometrically finite or Fuchsian (e.g., 
see \cite{Sullivan79}, \cite{Sullivan84}).
The direction $\dim(\Lambda_c(G)) \leq \delta(G)$ is easy and 
well known and the opposite inequality is fairly simply (and 
known) when the
Poincar{\'e} series diverges at the critical exponent.
We should also note that \thmref{delta=conical1} holds for 
a discrete group of M{\"o}bius transformations acting on the 
hyperbolic ball in any dimension as well as the rank 1 case in 
general.
 
It is known that the conical limit set corresponds to 
the geodesics starting at $0 \in \Bbb B$ which return to
some compact subset of $M$ infinitely often. Our proof 
shows that for any $\epsilon>0$ there is a subset of 
$\Lambda_c(G)$ of dimension $\geq \delta(G) - \epsilon$
which corresponds to geodesics which never leave some 
some compact subset of $M$. Thus our proof shows that 
if $M$ is a hyperbolic mainifold with finitely generated
fundamental group then 
for any  $x\in M$, the set of directions 
corresponding to geodesics rays starting at $x$ which have 
compact closure has dimension $\delta(G)$.

\thmref{thm1}  uses nothing about Kleinian groups except 
the definitions  but does make use of some nontrivial results
from the theory of conformal mappings and rectifiable sets. 
We shall actually give a long list  of conditions which are 
equivalent to $\Omega$ not being a disk (see the end of \secref{parabolics}.
Among them:
\begin{enumerate}
\item $\delta(G) >1$.
\item $\partial \Omega$ has tangents almost nowhere (with respect 
to harmonic measure).
\item $\partial \Omega$ fails to have a tangent somewhere.
\end{enumerate}
The elementary groups have to be excluded in \thmref{thm1} and 
\thmref{delta=dim1} because a cyclic group consisting of 
parabolics has a one point limit set, but $\delta (G) =1/2$.

\thmref{thm2} is essentially a corollary of \thmref{thm1} but 
requires the Ahlfors finiteness theorem and the Klein-Maskit 
combination theorems to reduce to the case considered in 
\thmref{thm1}.
The result was first formulated by Bowen in \cite{Bowen79}
in the case of quasi-Fuchsian groups with no parabolics. The convex, 
co-compact Kleinian  case is proven  in \cite{Sullivan79} and 
\cite{Braam88}. See also \cite{Sullivan82}.
The general geometrically finite case is 
proven  in \cite{Canary-Taylor93}. In this 
paper we complete the discussion by including 
the geometrically infinite groups.

We should also note that \thmref{thm1} and \thmref{thm2}
could be deduced from \thmref{gi->dim=2} and the known 
results for geometrically finite groups. However, our 
proofs of these results  seem new even in the geometrically finite case, 
so we have  included them in this paper.

\thmref{gi->dim=2} is the first place we have to invoke some 
nontrivial results about $3$-manifolds.
We need a result 
of Sullivan \cite{Sullivan87} that $\lambda_0 = \delta(G)(2-\delta(G))$ where
$\lambda_0$ is the lowest eigenvalue for the Laplacian on 
$M= {\Bbb B} / G$.
We  also need 
exponential decay (in time)  for the heat kernel of a
manifold with first eigenvalue bounded away from zero.

\thmref{gi->dim=2} was previously known in special cases.
Examples of groups with $\dim(\Lambda(G))=2$ were constructed by 
Sullivan in \cite{Sullivan81}, and Canary  \cite{Canary93}
proved \thmref{gi->dim=2}  holds
if $  M={\Bbb B}/G$ is a ``topologically tame'' manifold
and such that the thin parts have bounded type (in particular,
if its  injectivity radius is bounded away from zero).
Our result shows these hypotheses are unnecessary.

Sullivan \cite{Sullivan84} and Tukia \cite{Tukia84} independently
showed that if $G$ is a geometrically finite group then 
$\dim(\Lambda(G)) <2$. Thus \thmref{gi->dim=2} 
implies that a finitely generated group is geometrically finite 
iff $\dim(\Lambda(G)) <2$.

If $M = {\Bbb B} / G$, the convex core of $M$ is defined to 
be $C(M) = C(\Lambda)/G$, where $C(\Lambda)$ is the convex hull 
of the limit set in $\Bbb B$.
Burger and Canary  prove in 
\cite{Burger-Canary93} that if $G$ is a geoemtrically finite,
$M = {\Bbb B} /G$ has infinite volume and $\delta(G) > 1$, then 
$$ 2- {2 \text{area}(\partial C(M)) \over \text{vol} C(M) }
   \leq \dim(\Lambda(G))
   \leq 2 - { K \over \text{vol}(C_1(M))},$$
where $C_1(M)$ is a unit neighborhood of $C(M)$. \thmref{gi->dim=2}
implies their result still holds in the geometrically infinite
case (finitely generated groups are geometrically 
infinite iff  $\text{vol}(C(M)) = \infty$; see \secref{convex-hull}).

To prove \thmref{Hdim=lsc} we use our earlier results and 
an estimate of Canary which bounds the first eigenvalue of $M$
in terms of the volume of the convex core  $C(M)$.
We also need the use the Margulis lemma.
These results are  only needed in the
case when $\Lambda(G)$ has positive area; otherwise \thmref{Hdim=lsc}
follows immediately from earlier results.
We will also see that the assumption that $G$ and the $\{G_n\}$ 
are Kleinian is unnecessary; the result holds for any finitely 
generated M{\"o}bius groups.

  \thmref{H=M1} and \thmref{delta=dim1} are
both corollaries of a result relating $\delta(G)$ to Minkowski
(or ``box counting'') dimension which we will state in 
\secref{Minkowski}. This result uses only \thmref{delta=conical1}, 
the Ahlfors finiteness theorem and some geometric arguments.
\thmref{H=M1} was proven in the geometrically finite case
by Stratmann and Urbanski \cite{Stratmann-Urbanski}, and our 
proof is quite simple in this case.  The geometrically infinite case
follows immediately from \thmref{gi->dim=2}, although there is 
also a proof which does not require \thmref{gi->dim=2}.

The Ahlfors conjecture states that the limit
set of a finitely generated discrete group of M{\"o}bius
transformations is either the 
whole sphere or has zero area. We do not address the 
Ahlfors conjecture in this paper, but if it is 
true several of the arguments given here would simplify. As
part of the proof of \thmref{gi->dim=2} we show that 
a finitely generated, geometrically infinite group with 
$\delta(G) <2$ must have a limit set with positive area.
Thus the Ahlfors conjecture implies $\delta(G)=2$ 
for any geometrically infinite
group. We do not know  an argument for the converse
direction, but both results are currently
known to be true for topologically tame groups, \cite{Canary93}.

In most of our results the case when $G$ has no rank 1  parabolics is 
easier because then the action of $G$ on $\Omega(G)$ has a 
compact fundamental polygon. When the parabolics introduce 
extra difficulties we usually give the proof first in the compact
case and then go on to the general case.
We will introduce the notion of ``good'' and ``bad'' horoballs.
The good horoballs can be treated very similarly to the compact case
but the bad horoballs require a little extra work (and can only 
occur in the geometrically infinite case).
 The rank 2 parabolics 
cause no difficulties since they do not correspond to 
cusps on $\Omega(G)/G$.

It might be worth noting that although the theory of hyperbolic 
$3$-manifolds is intricate and highly developed,  we 
we only need a few simple facts which we prove here. The only 
results we really use without proof are the Margulis lemma, Davies'
Gaussian upper bound on the heat kernel 
and Sullivan's theorem  relating
the lowest eigenvalue for the Laplacian on $M ={\Bbb B} /G$ to 
the Poincar{\'e} exponent $\delta(G)$.

The remaining  sections of this paper are organized as follows:
\begin{description}
\item[\secref{sec-conical}] 
We prove \thmref{delta=conical1} and deduce from the proof that 
$\delta(G)$ is lower semi-continuous with respect to algebraic 
convergence.
\item[\secref{beta}]
This is first of several sections devoted to \thmref{thm1}
and \thmref{thm2}.
We define the ``$\beta$'''s which measure 
the distance of a set to a line and 
we prove that a uniformly wiggley 
set (i.e., a set where the $\beta$'s are bounded away 
from zero) has dimension larger than one. 
\item [\secref{big-Sch}] We recall that large Schwarzian 
for a conformal map implies large $\beta$'s for the image 
domain. This implies that if $G$ has no parabolics then $\Lambda(G)$
is either a circle or
is uniformly wiggley. Thus we obtain the compact case of 
\thmref{thm1} as a corollary.
\item[\secref{non-diff}] We prove  that if a simply connected invariant
component is not a disk then the boundary has tangents almost 
nowhere.  Here we need some 
known  results about conformal mappings and Schwarzian derivatives.
We also prove a result on the length of level lines that is needed
in the next section.
\item[\secref{parabolics}] We prove  \thmref{thm1} with parabolics.
\item[\secref{round}] We show that if $\Omega(G)$ is a union of 
infinitely many disks, then $\dim(\Lambda(G)) >1$. This is easy and previously 
known.
\item[\secref{degen}] We show that the limit sets of degenerate
groups (groups where $\Omega(G)$ is connected and simply connected)
are uniformly wiggley. This gives an alternate proof of 
\thmref{thm1} in this case.
\item[\secref{pf-thm2}] We prove \thmref{thm2} and deduce a 
variety of corollaries. 
\item[\secref{G&B-horoballs}]
We prove $\Lambda(G)$ is uniformly perfect and  introduce the
 idea of ``good'' and ``bad''
horoballs and prove some simple facts about them.
\item[\secref{convex-hull}] Here we gather together the 
facts about the convex hull of a hyperbolic $3$-manifold
that we will need in later sections.
\item[\secref{geo-inf}] We prove \thmref{gi->dim=2}.
\item[\secref{LSC-dim}] We prove \thmref{Hdim=lsc}.
\item[\secref{UMdim=BT}] We introduce Minkowski and upper 
Minkowski dimension and relate it to the Besicovitch-Taylor
index. 
\item[\secref{Minkowski}] We prove \thmref{H=M1} and \thmref{delta=dim1}
by proving the  following result: if $G$ is finitely generated 
and  $\Lambda(G)$ has zero area then 
$\delta(G)$ equals the upper Minkowski dimension of $\Lambda(G)$.
\item[\secref{Teich}] We deduce some corollaries of our 
results in the special case when the groups belong to 
$\overline{T(S)}$, the 
closure of the Teichm{\"u}ller space of a finite type surface $S$.
For example, $\dim(\Lambda(G))$ is a 
lower semi-continuous function on  $\overline{T(S)}$
and is continuous everywhere except at the 
geometrically finite cusps, where it is discontinuous.
\end{description}

The first author thanks many people for helpful 
conversations and suggestions; Michael Anderson,  
Itai Benjamini, Mladen Bestvina, Dick Canary, Irwin Kra,  Bernie Maskit,
Yair Minsky and Yuval Peres among others. 
He  particularly thanks Ed Taylor for 
explaining the results of \cite{Canary-Taylor93} and 
for suggesting the problem of proving the  Hausdorff dimension
of degenerate limit sets is strictly greater than 1. It was 
by considering this problem that we were led to the other
results in this paper.

\def\uMdim{\overline{\text{Mdim}}}
\def\lMdim{\underline{\text{Mdim}}}
\def\Mdim{\text{Mdim}}
\def\Ncov{N_{\text{cov}}}
\def\Cal{\cal}

\def\disk{\Bbb D}
\def\sphere{\overline {\Bbb C}}
\def\dist{\text{dist}}
\def\diam{\text{diam}}
\def\dim{\text{dim}}
\def\Dbar{\overline{\partial}}

\section{ The conical limit set}\label{sec-conical}
First we recall the definition of Hausdorff dimension. 
Given an increasing function $\varphi$ on $[0, \infty)$, 
we define  
$$ H_\varphi^\delta (E)=   \inf\{\sum \varphi(r_j):
E \subset \cup_j D(x_j, r_j), r_j \leq \delta\},$$
and 
$$ H_\varphi (E)= \lim_{\delta \to 0} H_\varphi^\delta (E).$$
This is the Hausdorff measure associated to $\varphi$. 
$H_\varphi^\infty$ is called the Hausdorff content. It is 
not a measure, but has the same null sets as $H^\alpha$.
When $\varphi(t) = t ^\alpha$ we denote the measure $H_\varphi$  by 
$H_\alpha$ and we define
$$ \dim(E) = \inf\{ \alpha: H_\alpha(E) =0 \}.$$
For $\alpha =1$ we sometimes denote $H_1$ by $\ell$ 
(for ``length'').
An upper bound for $\dim(E)$ can usually be produced 
by finding appropriate coverings of the set. We will be 
more interested in finding lower bounds. The usual 
idea is to construct a positive measure $\mu$ on $E$
which satisfies $\mu(D(x,r)) \leq C r^\alpha$. This implies
$\dim(E) \geq \alpha$ since for any covering of $E$ we have
$$ \sum_j r_j^\alpha \geq C^{-1}\mu(D(x_j, r_j) \geq
C^{-1} \mu(E) > 0.$$

Next we recall the definition of the critical exponent of Poincar{\'e}
series of $G$, 
$$\delta(G) = \inf\{ \alpha: \sum_{g \in G}(1-|g(0)|)^ \alpha < \infty\}.$$
It is easy and well known that 
$$ \dim(\Lambda_c) \leq \delta(G),$$
so we will only prove the other direction.
If the Poincar{\'e} series diverges at the critical exponent
$\delta(G)$, then the Patterson measure  (\cite{Patterson76}, 
\cite{Patterson84}, \cite{Nicholls88}, \cite{Sullivan84})
lives on the conical limit set $\Lambda_c(G)$ and satisfies the 
estimate
$$ \mu(B(x,r)) \leq C r^{\delta(G))}.$$
So the theorem is known in this case and we must only treat the 
case where the Poincar{\'e} series converges at the 
critical exponent.
In this case the Patterson measure  is supported on $\Lambda(G) 
\setminus \Lambda_c(G)$. Our proof does not use Patterson 
measure or any property of divergence or convergence 
at the critical index. It does construct a 
measure on the conical limit set which mimics the behavior of the 
Patterson measure.

In this section we will prove,

\begin{thm} \label{delta=conical}
Suppose $G$ is a group of M{\"o}bius transformations with more
than two limit points.
Then
$$ \delta(G)  \leq \dim(\Lambda_c).$$
\end{thm}

\begin{pf}
Let $\{ z_n\}$ denote the orbit of $0$ under $G$ in the 
hyperbolic $3$-ball, $\Bbb B$. Let $ \delta = \delta(G)$
be the 
critical exponent for the Poincar{\'e} series of $G$ and
let $\epsilon >0$.  Choose a point $x \in \Cbar = 
\partial \Bbb B$ so that
$$ \sum_{j:|z_j-x|<r} (1-|z_j|)^{\delta -\epsilon} = \infty ,$$
for every $r >0$ (here $|z-x|$ denotes the spherical metric).
 We can do this by a simple compactness argument.
Since $G$ is non-elementary, $x$ is not fixed by every 
element of $G$. Therefore we can choose elements $\{g_1, 
\dots, g_4\} \in G$ so that  $x_i = g_i(x)$ are all 
distinct. Fix $r>0$ to be so small that 
the balls $B_i$ on $\Cbar$ (in the spherical metric) of radius $2r$
around the points $\{x_1, \dots x_4\}$ are pairwise disjoint.

Suppose $M, N$ are large numbers (to be chosen below depending 
only on $G$ and  $r$).
Let $A_n= \{ z \in {\Bbb B}: 2^{-n-1} \leq 1-|z| < 2^{-n} \}$.
 If it were true that 
$$ \sum_{j: z_j \in B_i  \cap A_n } (1-|z_j|)^{\delta
- 2\epsilon} \leq M,$$
for all large enough $n$, then 
\begin{eqnarray*}
\sum_{j: z_j \in B_i }(1-|z_j|)^{\delta - \epsilon} 
&\leq& 
C \sum_n 2^{-n \epsilon} \sum_{j: z_j \in B_i \cap A_n}
 (1-|z_j|)^{\delta - 2\epsilon}\\
& \leq&  C M \sum_n 2^{-n \epsilon}\\
&<& \infty.
\end{eqnarray*}
This is a contradiction, so we must have 
$$ \sum_{j: z_j \in B_i \cap A_n} (1-|z_j|)^{\delta
- 2\epsilon} \geq M,$$
for infinitely many values of $n$. For each $i = 1,2,3,4$, fix a value 
of $n_i$  (depending on $M$ and hence of $r$) 
for which this inequality holds. 

Since the $z_j$'s make up the orbit of a single point, they are
uniformly separated in the hyperbolic metric of $\Bbb B$. Thus for 
any $A< \infty$ we may split the sequence into a finite number $B$ of 
sequences (depending on $A$) each of which is separated by at least 
$A$ in the hyperbolic metric. Therefore, to each point $x_i$ we may 
associate a collection of points $\Cal G_i(0) \subset \{z_j\}$ such that
$$  \Cal G_i(0) \subset B_i \cap A_{n_i},$$
$$ z,w \in  \Cal G_i(0) \text{ implies }  |z-w| \geq 3 N 2^{-n_i},$$
$$ \sum_{j: z_j \in  \Cal G_i(0)} (1-|z_j|)^{\delta
- 2\epsilon} \geq M/B.$$

For each $z \in  \Cal G_i(0)$ let $z^*$ denote its radial projection onto
the sphere $\Cbar = \partial \Bbb B$. For $z \in \Bbb B$, let 
$$ B(z_j) = B(z_j^*, N (1-|z|)).$$
By hypothesis, the balls $B(z_j)$ are disjoint for all 
$ z_j \in \cup_{i=1}^4  \Cal G_i(0)$.

Since the balls $\{ B_1, \dots, B_4\}$ have disjoint
doubles, any sufficiently small disk (depending on $ r$)
 can intersect at most one of the balls.
For any point $z = g(0)$ in the orbit of $0$  consider the 
four balls $\{ g(B_1), \dots, g(B_4)\}$. The preceding 
statement implies that if $N$ is sufficiently large (depending 
only on $r$) then at most one of these balls
can intersect ${\Bbb B} \setminus B(z)$. This determines our 
choice of $N$. Therefore, at least three
of the balls are contained in $B(z)$. Without loss of generality, 
assume they are $g(B_1), g(B_2)$ and $g(B_3)$.

The M{\"o}bius transformation $g$ has bounded distortion as 
a map from $S^2$ to itself except possibly
at one point. More precisely, 

\begin{lemma}\label{ex-point}
Suppose $r>0$. There is a $C< \infty$ (depending only on $r$)
such that given any M{\"o}bius transformation $g$ of 
$S^2= \Cbar$ to itself we have 
$$  C^{-1} (1-|g(0)|) \leq |g'| \leq C (1-|g(0)|),$$
except possibly on a disk $D$ of radius  $r$ (both the derivatives
and the disk are taken with respect to the spherical metric).
\end{lemma}
 
\begin{pf}
We may assume $g(0)\ne 0$ since otherwise the lemma is easy.
Let $z$ be the radial projection of $g(0)$ onto $S^2$ and 
choose $R$ so big (depending only on $r$) so that 
$$\omega(g(0), S^2 \setminus B(z,R(1-|g(0)|), {\Bbb B}) \leq r^2.$$
Then $D= g^{-1}(B(z,R(1-|g(0)|))$ is a disk of radius less than $r$ and 
$|g'|$ is comparable to $1-|g(0)|$ on its complement.
\end{pf}

So for any $g \in G$, at least two of the disks (say $B_1$ and $B_2$) are
bounded away from this point  so we get 
$$  C^{-1} (1-|g(0)|) \leq |g'| \leq C (1-|g(0)|),$$
on $B_1$ and $B_2$ with constants depending only on  $r$.

Note that if we choose $n_i$ large enough (depending only on $N$) 
we may assume 
$$ {1-|z| \over 2 N} \geq  1-|z_j| \geq {1-|z| \over C N},$$
for some uniform $C$ depending only on $G$ and $r$.

Now for the orbit point $z= g(0)$ define $\Cal G(z) = g(\Cal G_1(0)) $.
Thus 
\begin{eqnarray*}
\sum_{z_j \in  \Cal G(z) } (1-|z_j|)^{\delta - 2 \epsilon}
&\geq& C^{-\delta + 2 \epsilon} (1-|z|)^{\delta - 2 \epsilon}
    \sum_{z_k \in  \Cal G_1(0)} (1-|z_k|)^{\delta - 2 \epsilon} \\
&\geq & C^{-\delta+2 \epsilon}  (1-|z|)^{\delta - 2 \epsilon} \frac MB \\
&\geq & C^{-2} \frac MB (1-|z|)^{\delta - 2 \epsilon} \\
&\geq & (1-|z|)^{\delta - 2 \epsilon},
\end{eqnarray*}
 where the last line holds if $M$ is  large enough. Since 
$C_1$  depends only on $r$ and $B$ depends only on  group $G$
(more precisely it only depends on the injectivity radius of 
$G$ at $0$), this determines our choice of $M$.

We have now constructed a set of points ${z_j} = \cal G(z)$
which  satisfy the following conditions: 
$$ z_j \in B(z^*, N(1-|z|)),$$
$$  B(z_j^*, 2N(1-|z_j|))
\cap B(z_k^*, 2N(1-|z_k|)) = \emptyset \text{ for } j \ne k,$$
$$\sum_j (1-|z_j|)^{\delta- 2 \epsilon} \geq (1-|z|)^{\delta - 2 \epsilon},$$
$$ {1-|z| \over 2 N} \geq  1-|z_j| \geq {1-|z| \over C N},$$
for some uniform $C$ depending only on $G$ and $r$ (because 
the points in $G_1(0)$ do and the map $g$ has uniformly bounded 
distortion on $G_1(0)$).	
It is now a standard argument to show that $\dim(\Lambda_c) 
\geq \delta - 2 \epsilon$. Briefly it goes as follows.

Construct generations of points starting with $\cal G_0 = \{ 0\}$, 
and for each $z \in \cal  G_n$, define points $\{ z_j\}$ in 
$\cal G_{n+1}$ as above. To each point $z \in \cal G =\cup_n \cal G_n$, 
associate the disk $ B_z=B(z^*, 2N(1-|z|))$. Then let
$$ E_n = \cup_{z \in \cal G_n} B_z,$$
$$ E = \cap_{n} E_n.$$
It is easy to see that $E \subset \Lambda_c(G)$.

Define a  probability measure $\mu$ on $E$ by setting $\mu(E_0)=1$,
and for $z \in \cal G_n$ with ``parent'' $z'\in \cal G_{n-1}$, set
$$ \mu(B_z) = {(1-|z|)^{\delta-2\epsilon} \over \sum_{w \in \cal G(z')}
(1-|w|)^{\delta-2\epsilon} } \mu(B_{z'}).$$
It is easy to see by induction that 
$$\mu(B_z) \leq (1-|z|)^{\delta- 2 \epsilon} \leq C \diam(B_z)^{\delta
-2 \epsilon},$$
for each $z$ in $\cal G$. We want to show this inequality is true for
any disk $D$. Let $D$ be any disk and let $D_0 = B_z$ be the 
lowest generation disk in our construction so that 
$D_0 \cap D \ne \emptyset$ but $D_0 \subset 2D$. Let $D_1$
be the parent of $D_0$. By the maximality of $D_0$ we have 
$D \subset 2 D_1$. Since $2 D_1$ is disjoint from any other 
balls of the same generation, 
\begin{eqnarray*}
\mu(D) &\leq & \mu(D_1)\\
& \leq& C \diam(D_1)^{\delta- 2\epsilon}\\
&\leq &  C (N C)^{\delta - 2 \epsilon}\diam(D_0)^{\delta- 2\epsilon}\\
&\leq &   2C (N C)^{\delta - 2 \epsilon} \diam(D)^{\delta- 2\epsilon}.
\end{eqnarray*}
This is the desired inequality (the constant in front is larger, but 
is uniform over all disks; the power is the same).

Now we simply note that
if $\{ D_j\}$ is any covering of $E$ by disks, then 
$$ 0< \mu(E) \geq \sum_j \mu(D_j) \leq C \sum \diam(D_j)^{\delta- 2\epsilon}.$$
By the definition of Hausdorff dimension, this implies 
$$ \dim(\Lambda_c(G)) \geq \dim(E) \geq \delta- 2 \epsilon.$$
Since $\epsilon$ was arbitrary, we get \thmref{delta=conical}.
\end{pf}

The proof of \thmref{delta=conical} has the following 
corollaries.

\begin{cor}
Suppose $r>0$ and let $C = C(r)$ and $N = N(r)$ be the constants
in \lemref{ex-point} and the proof above.
Suppose $G$ is a group of M{\"o}bius transformations on 
$\Bbb B$ and suppose there are three disjoint balls $B_1, B_2, B_3$
of (spherical) radius $r$ and a collection of points 
$F\subset G(0) \cap A_n$ which
satisfy
$$ z,w \in F \text { implies } |z-w| \geq  N2^{-n},$$
$$ \sum_{z\in F \cap B_i} (1-|z|)^\alpha \geq C^{-2}.$$
Then 
$$ \delta(G) \geq \alpha.$$
\end{cor}

Once we have the conditions in the hypothesis, the proof of 
\thmref{delta=conical} proves the lemma. If $G$ is a group
then it satisfies these conditions for every $\alpha < \delta(G)$.
If $\tilde G$ is another group which is very close to $G$
(say $G$ and $\tilde G$ have generators which are very 
close in $\text{PSL}(2, \Bbb C)$) then $\tilde G$ will 
also satisfy these conditions (since they only involve a 
finite number of elements in the group). Thus

\begin{cor}
Suppose $G$ is a M{\"o}bius group generated by $\{g_1, \dots, g_n\}$. 
Given any $\delta_0 >0$ there is an $\epsilon_0 >0$ (depending 
only on $\delta$ and $G$) such that if $\tilde G$ is a group
containing elements $\{\tilde g_1, \dots, \tilde g_n\}$ with 
$\|g_i - \tilde g_i\|< \epsilon$ (as elements of 
$\text{PSL}(2, \Bbb C)$) then 
$$ \delta(G') \geq \delta(G) - \delta_0.$$
\end{cor}

Suppose $\{G_n\}$ is a sequence of $m$-generated  M{\"o}bius groups each with 
a specific listing of its generators $G_n=\{ g_{1 n}, \dots g_{m n}\}$.
  We say 
that $G_n$ converges {\it algebraically} to a Kleinian group 
$G $ with generators  $\{g_1, \dots, g_m\}$ 
if $ g_{jn} \to g_j$ for each $1\leq j\leq m$, as
elements of $\text{PSL}(2, \Bbb C)$. See \cite{Jorgensen-Marden90}.
If we identify groups with points in $\text{PSL}(2, {\Bbb C})^m$, this
is just convergence in the product topology.

\begin{cor} \label{delta=lsc}
If $\{ G_n\}$ is a sequence of discrete M{\"o}bius groups 
 converging algebraically to $G$, then 
$$ \liminf_{n \to \infty} \delta(G_n) \geq \delta (G).$$
\end{cor}
 
This says that $\delta(G)$ is lower semi-continuous 
with respect to algebraic convergence.
Strict inequality is possible
even for sequences of Kleinian groups
(e.g., one can choose a sequence $\{G_n\}$ of geometrically finite groups 
in $T(S)$ converging to a a geometrically finite cusp group $G$
so that $(\delta(G_n) \to 2$ but $\delta(G)<2$).

We should also note that \corref{delta=lsc} is still true
even if the groups involved are not discrete (assuming that
we define $\delta(G)$ appropriately). If $G$ is a finitely 
generated M{\"o}bius group which is not discrete, then its 
closure is a closed subgroup of $\text{PSL}(2, \Bbb C)$ and 
hence a Lie subgroup. The only possible limit sets  (i.e., 
accumulation sets on $S^2$ for the orbit of $0$) are: zero,
one or two points, a circle or the whole sphere. Let 
us set $\delta(G)$ to be $0,1$ or $2$ in these three cases 
respectively. In each case, it is easy to check that 
\corref{delta=lsc} still holds if we allow either the 
$\{G_n\}$ or $G$ to be non-discrete groups.

\section{Big $\beta$'s imply big dimension} \label{beta}

Although Kleinian groups are closely associated with 
3-manifolds, our approach to \thmref{thm1}
is completely two dimensional, 
exploiting connections between conformal mappings, 
the Schwarzian derivative
and certain geometric square functions which first arose
in the study of the Cauchy integral.
 A dyadic square 
$Q$ in the plane is one of the form 
$ Q=[2^{-n}j, 2^{-n}(j+1)] \times [2^{-n}k, 2^{-n}(k+1)]$.
The side length of a square will be denoted $\ell(Q)$. 
For a positive number $\lambda >0$, we let $\lambda Q$ 
denote the square concentric with $Q$ but with side length
$\lambda \ell(Q)$, e.g., $2Q$ is the ``double'' of $Q$.

 Given a set 
$E$ in the plane  and  a square $Q$ we define 
$\beta(Q)$ and $\delta (Q)$ as 
$$ \beta(Q) = t^{-1} \inf_{ L \in \cal L} \sup_{z \in E \cap 3Q}
\text{dist}(z, L),$$
where $\cal L$ is the set of all lines  $L$  intersecting $Q$ and 
$$ \delta(Q) = \ell(Q)^{-1} \inf_{ L \in \cal L} H(E\cap 3Q, L \cap 
3Q),$$
where $H(E,F)$ denotes the Hausdorff distance between the sets,
$$H(E,F) = \max_{x \in E} \text{dist}(x,F)+
\max_{y \in F} \text{dist}(y,E).$$

The second author proved in \cite{Jones90} that  a planar set $E$ lies
on a rectifiable curve iff 
$$ \sum_Q \beta(Q)^2 \ell (Q) < \infty,$$
where the sum is over all dyadic squares in the plane.
 We shall show 

\begin{thm} \label{big-beta}
Suppose $E$ is a closed, connected set in the plane and 
$\beta_E(Q) \geq \beta_0 >0$ for every square $Q$ with 
$Q \cap E \ne \emptyset$ and $\ell(Q) \leq \text{diam}(E)$. 
Then $\text{dim}(E) \geq 1 + C\beta_0^2$, where $C$
is an absolute constant.
\end{thm}

This is sharp (except for the choice of $C$) as can be 
seen by considering the standard iterative construction of a snowflake
curve where each segment of length $r$ is divided into four 
equal segments and the middle two are  replaced by segments
of length $ r(\frac 14 +\beta^2)$. 

We will see later that if $G$ has no parabolics then $\Lambda(G)$
has this property (i.e., is uniformly wiggley'').
However, a general limit set
need not have $\beta$'s bounded away from zero everywhere (this 
fails near doubly cusped parabolic points), but they are 
large often enough to push the dimension above $1$.


{\it Proof of \thmref{big-beta}.}
Suppose $E_0$ is a connected compact set with large 
$\beta$'s, $Q$ is a square of side length $\ell(Q) 
=r = 2^{-N}$  with 
$\frac 13 Q \cap E \ne \emptyset$ (here and later
$\lambda Q$  denotes the square concentric with 
$Q$ of side length $\lambda \ell(Q)$).
Our first objective is to show that for small $\epsilon >0$ we can 
find more than $1000 \epsilon^{-1}$ subsquares of $Q$ 
of sidelength  $\epsilon \ell (Q)$ with 
disjoint doubles and such that each contains a point of 
$E_0$ in its middle third. We will then apply the same argument
to each subsquare and use the resulting nested collection to 
squares to build a Frostman measure on $E_0$.

 Define $E = (E_0 \cap Q) \cup \partial Q$. Note that $E$ is 
connected. Fix some integer $n$ so that $2^{-n} < r$ (possibly
much smaller). Because $E$ is connected it is easy to check that 
there are more than $\frac 13 r 2^{n}$ dyadic subsquares $\{Q_j\}$ of $Q$
which lie in $\frac 23 Q$ and such that $Q_j \cap K \ne
\emptyset$ (e.g., consider concentric ``annuli'' between 
$\frac 13 Q$ and $\frac 23 Q$ made of squares of size 
$2^{-n}r $; there are $\frac 13 r 2^{n}$ such and each must 
intersect $E$).

From this we deduce
$$ \sum_{Q: \ell(Q) = 2^{-n}} \beta_E^2(Q) \ell (Q) \geq \frac 13
 \beta_0^2 r.$$
Thus for any integer $k \geq 1$ (recall $r = 2^{-N}$),
$$ \sum_{n=N+1}^{N+k} \quad \sum_{Q: \ell(Q) = 2^{-n}} \beta_E^2(Q)
 \ell (Q) \geq \frac 13 k\beta_0^2 r.$$
Suppose $\Gamma_k$ is the shortest curve in the plane with the 
property that for each $z \in E$ we have $\text{dist}(z,\Gamma_k)
\leq 2^{-N-k}$. It is fairly easy to check that
that for squares with $\ell(Q) \geq 
10 \cdot 2^{-N-k} \beta_0^{-1}$, we have 
$$ \beta_{\Gamma_k}(Q) \geq \frac 12 \beta_E(Q) \geq \frac 12 \beta_0.$$
Therefore 
$$ \sum_{n=N+1}^{N+k} \quad \sum_{Q: \ell(Q) = 2^{-n}} \beta_{\Gamma_k}^2(Q)
 \ell (Q) \geq ( k- 10\beta^{-1})\beta_0^2 r.$$
Choose $k> 20 \beta^{-1}_0$.  Then the term on the 
right is $\geq \frac 16 k \beta_0^2 r$.

 By the second author's characterization of 
rectifiable curves in \cite{Jones90}, the length of 
$\Gamma_k$ is at least
$$ \ell(\Gamma_k) \geq C_0 \beta_0^2 k r,$$
for some absolute constant $C_0$. We claim that this 
implies that there are more than $(C_0 \beta_0^2 k r - 4r)C2^{N+k}$
boxes $\{Q_j\}$ of side length $2^{-N-k}$ such that
$\frac 13 Q_j \cap E \ne \emptyset$.

To prove this claim let $\{ z_j\}$ be a collection of points on
$\Gamma_k$ so that $j \ne k$ implies $|z_j-z_k| \geq 
2^{-N-k}$, but so that $\cup_j B(z_j, 2^{-N-k+2})$ covers
$E$. Let $\cal C$ be the collection  
$\{ \partial Q_k\}$ of all dyadic squares of size $2^{-N-k}$
contained in $Q$ which  contain some $z_j$. Let $\Gamma =
\cup_{Q_j \in \cal C} \partial Q_j \cup \cup_j S_j$ be the 
union of  the boundaries of these cubes, together with 
segments $S_j$ which connect $\partial Q_j$ with the point 
$z_j$. Then obviously
$$ \ell (\Gamma) \leq 6 \cdot 2^{-N-k} |\cal C|.$$
Since $\Gamma$ has the property that it passes 
within $2^{-N-k}$ of every point of $E$ and since $\Gamma_k$ 
was defined to be the shortest such curve, we
must have
$$  6 \cdot 2^{-N-k} |\cal C| \geq 
 \ell (\Gamma)  \geq \ell(\Gamma_k) \geq C_0 \beta_0^2 k 2^{-N}.$$
and hence
$$ |\cal C| \geq \frac 16 C_0 \beta_0^2 k 2^k.$$
So if  $k$ is large enough (e.g., $k > 60000 \beta_0^{-2} C_0^{-1}$)
Now set $\epsilon = 2^{-k}$. Then  $|\cal C| \geq 10000 \epsilon^{-1}$, 
and  consist of disjoint dyadic squares of size $\epsilon \ell(Q)$
each of which contains a point of $E$. By throwing away $9/10$'s
of the squares we can assume the remaining ones have disjoint
triples, and this is what we wished to prove.

To finish the proof of the theorem, we build nested generations
of squares using the construction above. The initial square 
$Q_0$ forms the first generation. The squares of size 
$\epsilon \ell(Q_0)$ constructed above form the first generation. 
In general, given a $n$th generation square containing a 
point of $E_0$ in its middle third, we construct $ 1000 
\epsilon^{-1}$  subsquares  as above (with disjoint triples 
and containing a point of $E_0$ in their middle thirds), and 
put these into the $(n+1)$st generation.

We then define a measure $\mu$ by assigning each $n$th generation
square equal mass (namely $(\epsilon/1000)^{n})$). Since a
$n$th generation square has size $\epsilon^n$, this measure 
satisfies 
$$ \mu(Q) \leq C \ell(Q)^{\alpha},$$
where 
$$ \alpha = {\log \epsilon - \log 1000 \over \log \epsilon}
= 1 + {\log 1000 \over \log \frac 1 \epsilon} >1,$$
and $Q$ is an $n$th generation square.
Since $\epsilon = 2^{-k}$ and $k \sim \beta_0^{-2}$  we 
get $ \log \epsilon^{-1} \sim \beta_0^{-2}$. This gives the 
estimate in the theorem. It only remains to
check that this inequality holds for all squares in the plane, 
but this is a standard argument. Thus
$$ \sum \text{diam}(D_j)^\alpha \geq \sum \mu(D_j) 
\geq \mu(E_0) > 0$$
for any covering $\{ D_j\}$ of $E_0$, and therefore
$\text{dim}(E_0) \geq \alpha >1$.$\qed$

\section{Large Schwarzian implies large $\beta$'s} \label{big-Sch}

The Schwarzian derivative of a locally univalent function $F$ is defined by
\begin{eqnarray*}
S(F)(z)  &=& [ {F''(z) \over F'(z)} ]' -{1\over 2}
[ {F''(z) \over F'(z)} ]^2\\
&=& [{F'''(z) \over F'(z)}] -{3\over 2} [ {F''(z) \over F'(z)} ]^2.
\end{eqnarray*}
If we write $F'= e^\varphi$ then it can be rewritten as
$$ S(F)(z) = \varphi'' - \frac 12 (\varphi ')^2.$$
Recall that $S(F) \equiv 0$ iff $F$ is a M{\"o}bius transformation
and that $S$ satisfies the composition law
$$ S(F \circ G) = S(F)(G')^2 + S(G).$$
In particular, if $G$ is M{\"o}bius then
$$ S(F \circ G) = S(F)(G')^2 $$
$$S(G \circ F) = S(F).$$
In addition, given an $\epsilon >0$, hyperbolic disk $D$ and a
compact neighborhood $K$ of $D$,
there is a $\delta >0$ so that $|S(F)|\leq \delta$  on $D$ implies
$F$ uniformly approximates a M{\"o}bius transformation on $K$ to
within $\epsilon$. 


We introduce some notation from the second author's paper
\cite {Jones89}.
For a domain $\Omega$ and a point $ z_0 \in \Omega$ let 
$ z_1 \in \partial \Omega$ be the closest point in the 
Euclidean metric and let $\theta_0 = \arg (z_0 -z)$.
 Let $\tilde z_1 = z_1 +  \delta (z_0, z_1)$ were 
$ 0< \delta < 1$ is fixed. Define 
$$ L^\delta_{z_0} = \{ z_1 + t \exp(i (\theta_0- \frac \pi 2)):
t \in {\Bbb R},  |t| \leq \delta^{-1} |z_1-z_0| \},$$
$$S^\delta_{z_0}  =\{ z: |z-z_1| \leq \delta^{-1} |z_1-z_0|, 
|\theta_0 - \arg(z- \tilde z_1)| \leq \frac \pi 2 \}.$$
Then $L^\delta_{z_0}$ is a sort of tangent line to 
$\partial \Omega$ (as seen from $z_0$) and $S^\delta_{z_0}$
is a half-disk which approximates $\Omega$ near $z_0$. We say 
that $\Omega$ satisfies condition $M(\delta)$  at $z_0$ if
either  there exists $z \in L^\delta_{z_0}$ such that 
$$ \dist(z, \partial \Omega) \geq \delta\dist(z_0, \partial 
\Omega),$$
or 
$$ S^\delta_{z_0} \cap \partial \Omega \ne \emptyset.$$
Thus if $M(\delta)$ is satisfied at $z_0$ then the 
boundary of $\Omega $ ``wiggles'' to order $\delta$
near $z_0$. The following is Theorem 2 of \cite{Jones89}

\begin{lemma}
Suppose $\Omega $ is simply connected and $ \Phi: \disk 
\to \Omega$ is a Riemann mapping.
For every $ \delta >0$ there is a $\epsilon >0$ so that if 
$s(z) = |S(\Phi)(z)|(1-|z|)^2 \geq \epsilon$ then 
$\Omega$ satisfies $M(\delta)$ at $\Phi(z)$.  Furthermore, 
$M(\delta)$ is satisfied at every point of $ \Omega$ iff
there are $C<\infty$ and $\epsilon>0$ so that $s(z) > \epsilon$ on 
a $C$-dense set of the disk (with respect to the hyperbolic metric).
\end{lemma}

The proof is a normal families argument. Since $\partial \Omega$
is connected it is easy to see that condition $M(\delta)$ 
implies that the square centered at $z_0$ of size 
$\delta^{-1} \dist(z_0, \partial \Omega)$ has 
$\beta_{\partial \Omega}(Q) \geq \delta^2$.  Thus large 
Schwarzian implies large $\beta$'s. Another  version (more 
precise than we need in this paper) is given by 
 
\begin{lemma} \cite{Bishop-Jones90} Suppose $E$ is compact, $\Omega 
= \sphere \backslash E$ and  $\Phi: \disk \to \Omega$ the 
universal covering map.  Suppose $\Phi$ is univalent on
$D(w,\epsilon (1-|w|))$.  Then
$$  |S(\Phi) (w)| (1-|w|^2 )^2\leq
C \epsilon^{-2}\sum_{n=0}^\infty \delta_E(2^n Q) 2^{-\mu n}  $$
where  $r = \dist(\Phi(w),\partial \Omega)$
The number $\mu$ satisfies $0 < \mu <1$ but can be taken as close to
1 as we wish.  The constant $C$ depends only on
the choice of $\mu$.
\end{lemma}

Using either of these estimates  we can prove the following 

\begin{cor} \label{Schwarz-dim}
Suppose $\Omega$ is simply connected and $\Phi: \disk \to 
\Omega $ is a Riemann mappings. Suppose also that there is 
a $C < \infty$ and $\epsilon >0$ so that every  point of 
the disk is at most hyperbolic distance $C$ from a point 
$z$ where $|S(\Phi)(z)(1-|z|^2)^2| \geq \epsilon$. Then there
is a $\beta = \beta(C, \epsilon)>0$ so that 
$\beta_{\partial \Omega}(Q) \geq \beta >0$ for every 
square $Q$ such that $\ell(Q) \leq \text{diam}(Q)$ and 
$ \frac 13 Q \cap \partial \Omega \ne \emptyset$. In particular, 
$\text{dim}(\partial \Omega) >1$ with estimates depending 
only on $C$ and $\epsilon$.
\end{cor}

\begin{pf}
We use the second lemma.
Let $E = \partial \Omega$ and let $Q$ be a square
  such that $\ell(Q) \leq \text{diam}(Q)$ and 
$ \frac 13 Q \cap E \ne \emptyset$. Let 
$\beta_0 = e^{-10 C}$. We may also suppose 
$\epsilon = 2^{- \mu N}$ for some integer $N$. If 
$\beta_E(Q) \geq  \epsilon \beta_0$, there is nothing to do, 
so assume $\beta_E(Q) \leq \epsilon \beta_0$.  By the Koebe
1/4 theorem, the  hyperbolic 
metric on $\Omega$ is comparable to $\text{dist}(z,E)^{-1} ds$, 
so by our hypothesis we can  find a point $z \in \frac 13 Q \cap \Omega$
with $  \epsilon \beta_0 \ell(Q) \leq\text{dist}(z,E)  \leq \epsilon
 \ell (Q)$, and such that if $ w = \Phi^{-1}(z)$ then 
$ |S(\Phi)(w)(1-|w|)^2| \geq \epsilon$.

Let $Q'$ be the square centered at $z$ with side length 
$2 \text{dist}(z, E)$. By the estimate above we know that
one of the squares $Q', 2Q',\dots 2^N Q'$ where $N  \sim
\log_2 \epsilon$ has $\delta_E \geq \epsilon$. Call 
this square $Q''$.
Since $z \in \frac 13 Q$ and $2^N \ell (Q) \leq \frac 14 \ell(Q)$,
 $Q''$ must be a subsquare of $Q$.   Since 
$\ell(Q'') \geq \beta_0 \epsilon \ell(Q)$, we deduce
$$\delta_E(Q) \geq \beta_0 \epsilon \delta_E(Q'')
\geq \beta_0 \epsilon^2 \equiv \beta.$$
This is what we wished to prove. The final statement is 
simply an application of \thmref{big-beta}. 
\end{pf}

Let $\Phi: {\Bbb D} \to \Omega$ be a Riemann mapping and 
let $\hat G = = \Phi \circ G \circ \Phi^{-1}$ denote 
the Fuchsian equivalent of $G$. We can now prove \thmref{thm1}
in the case when $\hat G$ has no parabolics.
In this case, the surface
$\Omega /G$ is compact so that any point of $\disk$
is within some bounded distance $C$ on any orbit of 
$\hat G$. If 
$\Omega$ is not a disk then $\Phi$ is not M{\"o}bius, 
so $S(\Phi)$ does not vanish identically. Choose a 
point $z$ so that $ s(z) =|S(\Phi)(z)(1-|z|^2)^2| = \epsilon \ne 0$.
Since $S(z)$ is constant on orbits of $\hat G$, we can 
apply \corref{Schwarz-dim} to deduce $\text{dim}(\partial 
\Omega) > 1$.

\section{Non-differentiability of the boundary} \label{non-diff}

In this section we describe our 
``non-differentiability'' result for  invariant components.
First we need some definitions which capture what we mean by 
``non-differentiable''.
A point $x \in \partial \Omega$ is called a (inner) tangent point 
of $\Omega$ if for any $\theta< \pi$, $x$ is the vertex of 
a cone in $\Omega$ with angle $\theta$, and this is not true 
for any $\theta > \pi$. The set of tangents is the same, up to 
a set of $1$-dimensional measure zero, as the set 
 $E\subset \partial \Omega$ of points which are
vertices of some cone in $\Omega(G)$.
 If $\Omega$ is simply connected and $x \in \partial 
\Omega$, then $x$ is called a twist point of $\Omega$ if
$$ \limsup_{z\to x, z \in \Omega} \arg(x-z)= + \infty
, \quad \liminf_{z\to x, z \in \Omega} \arg(x-z)=  -\infty.$$
The harmonic measure $\omega$ for the domain $\Omega$
is the push forward of Lebesgue measure on 
the circle under a Riemann mapping onto $\Omega$ (the measure 
depends on the choice of Riemann mapping, but its null sets
do not).  A probability measure $\mu$ is said to be singular 
with respect to  a measure $\nu$ if it gives full measure
to a null set for $\nu$.

\begin{thm} \label{cone}
Suppose $G$ is  a finitely generated Kleinian group with
simply connected invariant component $\Omega$. Then the 
following are all equivalent.
\begin{enumerate}
\item $\Omega$ is not a  disk.
\item The set of inner tangent points of $\partial \Omega$ has 
zero $1$-dimensional measure.
\item Harmonic measure for $\Omega$ is singular to 
$1$-dimensional measure.
\item Almost every (with respect to harmonic measure) point
of $\partial \Omega$ is a twist point.
\end{enumerate}
\end{thm}

Conditions (2), (3) and (4) are 
equivalent for any simply connected domain and they each 
imply (1). Thus the  main point is to show
(1) implies any of (2), (3) or (4).
  Using  known results 
on the behavior of harmonic measure (to be described below)
this result implies,

\begin{cor}
Suppose $G$ is finitely generated and $\Lambda(G)$ is not a
circle. Then the harmonic measures for different
components of $\Omega(G)$ are mutually singular.
\end{cor}

This last corollary is well known in many cases. For example, 
when $G$ is a quasi-Fuchsian group it corresponds exactly to 
Mostow's theorem: given two  finitely generated, first kind
Fuchsian groups on the unit 
disk and a homeomorphism $\phi$ of the boundary which conjugates 
the actions on the circle, then either  $\varphi$ is M{\"o}bius
or $\varphi$ is singular (i.e., maps full measure to zero measure).
A stronger version of Mostow's  result is proven in 
\cite{Bishop-Steger93} where it is shown that  such a
$\varphi$ is either M{\"o}bius or maps a set of dimension
$<1$ to the complement of a set of dimension $<1$.

If $\Omega$ is a component of $\Omega(G)$ then
the stabilizer of $\Omega$ in $G$ is also finitely 
generated (this is a consequence of the Ahlfors finiteness
theorem \cite{Ahlfors64}, \cite{Ahlfors69}).
 Thus if $\Omega$ is simply connected and  $\Phi: \disk \to 
\Omega_0$ is a Riemann mapping, $\hat G = \Phi \circ
G \circ \Phi^{-1}$ is a finitely generated Fuchsian
group of the first kind (i.e., its limit set is the 
whole circle) called the Fuchsian equivalent of $G$ on 
$\Omega$.
We can actually prove  \thmref{cone} using the weaker
hypothesis that the Fuchsian equivalent of $G$ on $\Omega$
has a non-tangentially dense orbit. This happens if it
is divergence type, i.e.,
$$ \sum_{\gamma \in \hat G} (1-|\gamma(0)|) = \infty.$$
This is the best one can expect since Astala and Zinsmeister
\cite{Astala-Zinsmeister90}
have shown that any convergence group
(i.e., $ \sum_{\gamma \in \hat G} (1-|\gamma(0|) < \infty$)
has a quasiconformal deformation to a
Kleinian group whose limit set is a rectifiable curve (not 
a circle).

Suppose $\Omega $ is simply connected and  $\Phi: {\Bbb D} \to \Omega$
is a Riemann mapping. By Plessner's theorem the unit circle 
$\Bbb T$ can be divided into two set $E_0, E_1$ such that almost 
every (Lebesgue measure) point of $E_0$, $\Phi'$ has a finite, non-zero 
non-tangential limit and almost everywhere on on $E_1$, 
$\Phi'$ is non-tangentially dense in the plane. Let 
$\omega$ denote the harmonic measure on $\Omega$, i.e., 
$\omega(E) = |\Phi^{-1}|/2 \pi$. By McMillan's theorem
\cite{McMillan69}, almost every point of $\partial \Omega$
is either the vertex of a cone in $\Omega$ or is a twist 
point (see introduction) and these two sets correspond a.e.
to $E_0$ and $E_1$ respectively via the map $\Phi$.  The set
of the cone points has $\sigma$-finite $1$-dimensional 
measure and on this set $\omega$ is mutually  absolutely 
continuous with $1$-dimensional Hausdorff measure.
Makarov  \cite{Makarov85} proved that there is a subset 
$F$ of the twist points which has zero $1$-dimensional measure
but the same harmonic measure as the full set of twist points.
Thus harmonic measure on a simply connected domain always 
``lives'' on a set of dimension $1$, regardless of the Hausdorff
 dimension of the  entire boundary. Pommerenke \cite{Pommerenke86}
proved that harmonic measure gives full measure to a set 
of $\sigma$-finite $1$-dimensional measure and Wolff \cite{Wolff93}
has even extended this to arbitrary  planar domains.

Given two disjoint simply connected domains $\Omega_1$, $\Omega_2$
with overlapping boundary $E = \partial \Omega_1 \cap 
\partial \Omega_2$, the two harmonic measures $\omega_1$, 
$\omega_2$ are mutually absolutely continuous on a subset 
$F \subset E$ iff almost every point (with respect to both 
measures) of $F$ is the vertex of two cones, one in each of 
the two domains \cite{BCGJ}. In this case, both harmonic measures are
also mutually absolutely continuous with respect to $1$-dimensional
Hausdorff measure on $F$. Thus the measures $\omega_1$
and $\omega_2$ are  mutually singular iff the set of 
``double cone points'' has zero $1$-dimensional measure.
In particular, if almost every point (with respect to $\omega_1$)
of $\partial \Omega_1$ is a twist point then $\omega_1$
is singular to $1$-dimensional measure and  singular to 
harmonic measure on any other disjoint domain. The converse is 
not true; there is a  closed Jordan curve so that the harmonic measures
for the two complementary domains are mutually singular 
but each is absolutely continuous with respect to 
$1$-dimensional measure \cite{Bishop87}. This is not possible 
for quasicircles however. A criterion  for singularity of 
harmonic measures of two general disjoint domains is given in 
\cite{Bishop91}.

There is also a connection between rectifiability and the 
Schwarzian derivative.

\begin{lemma} \cite{Bishop-Jones93} \label{Sch-length}
There is a $C < \infty$ such that if
 $\Phi: {\Bbb D} \to \Omega$ is univalent and $\partial \Omega$
is rectifiable (i.e., has finite $1$-dimensional measure) then 
$$ \iint_{{\Bbb D}} |\Phi'(z)| |S(\Phi)(z)|^2 (1-|z|^2)^3 dxdy 
< C \ell( \partial \Omega).$$
The same holds if $\Phi$ is only defined on a Lipschitz subdomain
of $\Bbb D$. If $\Omega $ is a quasicircle  and $\dist(\Phi(0), 
\partial \Omega) \sim \diam (\Omega)$ then the two sides are
comparable with a constant depending only on the quasicircle 
bound. 
\end{lemma}

This is Lemma 3.7 of \cite{Bishop-Jones93} and is an easy 
computation involving Green's theorem and basic estimates
for univalent mappings.

From the lemma we now  deduce \thmref{cone}. 
Suppose $\Omega$ is a simply connected component of $\Omega(G)$.
If $\Omega$ is not a circle, then $S(\Phi)$ is non-zero 
somewhere in the disk and therefore the invariant quantity 
$ s(z) = |S(\Phi)(z)| (1-|z|^2)^2$ equals $\epsilon >0$ on 
some orbit $\{ z_n\}$ of the Fuchsian equivalent $\hat G$.
Since $\hat G$ is a finitely generated Fuchsian group of the 
first kind, the orbit of any point is non-tangentially 
dense, i.e., for almost every point $x \in \Bbb T$ there
is a sequence of points in the orbit approaching $x$ through
a Stolz cone with vertex $x$. Alternatively, if we associate
to each orbit point $z_n$ the interval $I_n$ of length $1-|z_n|$
centered and $z_n /|z_n|$, then almost every point if the
circle is in infinitely many of the intervals $ \cal F=\{I_n\}$.

Now suppose $\Phi'$ has non-tangential limits on a set 
$E \subset \Bbb T$ of positive Lebesgue measure. Then by 
taking a union of small cones we could construct a ``sawtooth''
domain $W\subset \Bbb D$ with  $|E \cap \partial W| \geq
\frac 12 |E|$ and such that $M^{-1} \leq |\Phi'| \leq M$ on 
$W$.
Let $ F=E \cap \partial W$. Because of the non-tangential
density of the orbit we can use the Vitali covering lemma
(e.g. page 109 of \cite{Wheeden-Zygmund}) to obtain
a disjoint covering of almost every point of $F$
 by intervals of the for $\{I_n\}$.
In fact, by repeated applications of covering lemmas we can 
find infinitely many collections $\cal F_k = \{ I^k_j\} \subset
\cal F$ each of 
which is a disjoint covering of almost all of $F$, and 
so that no interval belongs to more than one collection.

For each orbit point $z_n$, let $D_n$ denote the disk centered
at $z_n$ of radius $\frac 12 (1-|z_n|)$. By the Koebe 1/4
theorem 
 \begin{eqnarray*}
 \iint_{D_n} |\Phi'(z)| |S(\Phi)(z)|^2 (1-|z|^2)^3 dxdy 
&\geq& C \epsilon^2 \iint_{D_n} |\Phi'(z)|(1-|z|)^{-1} dxdy \\
&\geq& C M ^{-1} \epsilon^2 (1-|z_n|)\\
&\geq & C M^{-1} \epsilon^2 |I_n|.
\end{eqnarray*}
Therefore 
 \begin{eqnarray*}
 \iint_{W} |\Phi'(z)| |S(\Phi)(z)|^2 (1-|z|^2)^3 dxdy 
&\geq & C M^{-1} \epsilon^2 \sum_{n: z_n \in W} |I_n| \\
&\geq & C M^{-1} \epsilon^2 \sum_k \sum_{I_n \in \cal F_k} |I_n| \\
&=& \infty.
\end{eqnarray*}
Therefore $\Phi(W)$ cannot have a finite length boundary by 
\lemref{Sch-length}. However, $|\Phi'|$ is bounded by $M$ 
on $W$, so $\Phi(W)$ must have a rectifiable boundary. 
This contradiction implies that $\Phi'$ cannot have 
non-tangential limits on any set of positive Lebesgue 
measure. Thus in terms of earlier remarks, almost 
every point of $\partial \Omega$ must be a twist point, 
the set of cone points has zero $1$-dimensional measure and 
harmonic measure on $\Omega$ is singular to harmonic measure 
for any simply connected domain which is disjoint from 
$\Omega$.

%

If we apply the proof not to the whole disk but to a 
``Carleson square'' $Q$ corresponding to an interval 
$I \subset \partial \Bbb D$, 
$$ Q = \{ re^{i \theta}: e^{i\theta} \in I , 1-|I| < r< 1\},$$
then  the proof shows that the length of 
$\Phi(rI)$ must tend to infinity. If we divide the circle 
into $N$ equal intervals, we can apply the proof to each 
interval and then by my taking the minimum growth rate, get
an estimate which is valid for all $N$ intervals. Since 
any interval of length $\geq 4\pi/N$ contains at least one
interval from our collection, the growth rate also holds 
for such an interval. Thus 

\begin{lemma}\label{long-length}
Suppose $G$ is finitely generated with invariant component
$\Omega$ which is not a  disk.
If $\Phi$ is a Riemann mapping onto  $\Omega$ and $I\subset
\partial \Bbb D$  then $\ell(\Phi(rI)) \to \infty $ as $r\to 1$
with estimates that only depend on $G$ and the length of $I$.
\end{lemma}

We finish this section with a remark concerning
Makarov's law of the iterated logarithm from \cite{Makarov85}.
Not only did Makarov prove that the harmonic measure for a simply connected
domain  gives full mass to a set of dimension 1 (as discussed
earlier), but also the following sharp converse.
If we define 
$$\varphi_C(t) = t \exp(- C \sqrt{ \log \frac 1t \log \log \log
\frac 1t }),$$
then there is an absolute constant $C_1$ so that harmonic 
measure  for any simply connected domain 
is absolutely continuous with respect to the 
Hausdorff measure for the gauge function $\varphi_{C_1}$. 
Conversely, there is a $C_2$ and examples of simply connected 
domains for which harmonic measure gives positive mass to 
a set of zero $H_{\varphi_{C_2}}$  measure. Such domains 
are called Makarov domains and have the worse possible 
behavior for harmonic measure.
In \cite{Jones89} the second author proved that $\Omega$ is 
a Makarov domain if its boundary  is ``wiggley'' on all
scales (in a sense made precise in that paper). It is 
easy to check that $\Omega$ satisfies his condition if 
every point of the disk is within a bounded hyperbolic distance
of a point where
$s(z)= |S(\Phi)(z)(1-|z|^2)^2$ is bounded away from zero.
We have already observed that if $\Omega$ is a component of 
a finitely generated group $G$ such that the Fuchsian 
equivalent $\hat G$ contains only hyperbolics, then this 
condition is satisfied. Thus all such domains are Makarov
domains. Presumably the same is true with parabolics, but
we have not attempted to prove this.

\section{ Proof of \thmref{thm1} with parabolics}\label{parabolics}

We have already proven \thmref{thm1} when the Fuchsian 
equivalent $\hat G$ of $G$ on $\Omega$ has no parabolic 
elements. In this section we will prove it assuming 
$\hat G$ does contain parabolics. This will also give a 
different proof for the co-compact case.

  This section is not 
strictly necessary for proof of \thmref{thm1}; we shall
later give two more proofs for the case of degenerate groups
and  \thmref{thm1} can be deduced
from this and the known results for geometrically finite 
groups.  However, this argument seems to be necessary 
to prove $\delta(G)>1$.
In fact, we will show $\delta(G) >1$ and use  the fact that 
if $\Omega$ is $G$-invariant, then $\partial \Omega = \Lambda(G)$.
Thus by \thmref{delta=conical1} we deduce \thmref{thm1}.

First we want to show  that it is enough to estimate the 
Poincar{\'e} series along an orbit in $\Omega(G)$ (rather 
than along the orbit of $0 \in \Bbb B$.
Recall that for $z \in \Omega(G)$ we set
$$d(z) = \dist(z, \partial \Omega(G)).$$

\begin{lemma} \label{B3=S2}
If $G$ is a Kleinian group, $0$ is the center of $\Bbb B$
 and $ z \in \Omega(G) \subset \partial \Bbb B$, then
$$d(g(z)) \sim 1-|g(0)|, \quad \text{for all } g \in G$$
with constants that depend on $z$ and $G$, but not on $g$.
\end{lemma}
 
\begin{pf}
Choose a ball $B \subset \Omega$ centered at $z_0$  so that
$\diam(B) \leq \frac 12 \dist (z_0, \partial \Omega$.
Let $ \omega_1 = \omega(0, B, \Bbb B)$ be the harmonic measure
of this ball in $\Bbb B$ with respect to the point zero and
let $ \omega_2$ be the harmonic measure
of $\frac 12 B$  with respect to $0$. Then by the conformal
invariance of harmonic measure,  for any $g \in G$,
$g(0)$ is the unique point $z$ so that
$$ \omega (z, g(B), {\Bbb B}) = \omega_1,$$
$$ \omega (z, g(\frac 12 B), {\Bbb B}) = \omega_2.$$
By our choice of $B$ and the Koebe $1/4$ theorem,
$g(\frac 12 B) \subset  \lambda g(B)$  for some $\lambda < 1$
independent of $g$. Therefore any  $z \in \Bbb B$ which satisfies the
two equalities above must satisfy
$$ |g(z_0) -z| \leq C\diam(g(B)),$$
$$ 1-|z| \geq \frac 1C \diam (g(B)).$$
This proves the lemma.
\end{pf}

Note that this implies that the accumulation set of any point
in $\Omega(G)$ is all of $\Lambda (G)$. In particular,
 if $\Omega$ is an  invariant component of $\Omega(G)$ then 
$\partial \Omega = \Lambda(G)$. 

Let $\hat G$ denote the Fuchsian equivalent of $G$ and 
suppose $0$ is a point where the Schwarzian derivative 
of the Riemann mapping ${\Bbb D} \to \Omega$ is large.
Suppose $z_j = g(0) \in \Bbb D$ is an orbit point of 
zero and let $I_j \subset \partial \Bbb D$ be the interval on 
the boundary centered at $z_j/ |z_j|$ with length  
$1- |z_j|$. We let $S_j$ denote the Carleson square with base 
$I_j$,
$$ S_j = \{ z: z/|z| \in I_j, 1-|I_j| \leq |z| <1\}.$$
Let 
$$ d_j = d(\Phi(z_j)) = \dist(\Phi(z_j), \partial \Omega)
\sim (1-|z_j|) |\Phi'(z_j)|.$$
 We will show 

\begin{lemma} \label{mainlem}
There is a $C < \infty$ (depending only on $G$)
 such that if   $g\in \hat G \setminus \{\text{Id} \}$,
and $z_0 =g(0)$ then there is a collection of orbit 
points $\{z_k\} = \{ g_k (0)\} \subset T_0$ such that
\begin{enumerate}
\item $ \sum_k d_k \geq  2  d_0,$
\item $(1-|z_k|) \geq (1-|z_0|)/C$.
\item the intervals $\{I_k\}$ are disjoint.
\end{enumerate}
\end{lemma}

Note that conditions (1) and (2) imply 
$$\sum_k d_k^{1+\epsilon} \geq  d_0^{1+\epsilon},$$
if $\epsilon$ is small enough (depending on $C$).
Using condition (3), we can break the orbit of $0 \in \Bbb D$ into 
generations $\cal G_n$ so that 
$$\sum_{z_k \in \cal G_n} d_k^{1+\epsilon} \geq  
\sum_{z_k \in \cal G_{n-1}} d_k^{1+\epsilon} \geq 1 ,$$
which implies
$$\sum_{z_k \in \hat G(0)} d_k^{1+\epsilon}  = \infty.$$
By \lemref{B3=S2}, this proves $\delta(G) >1$.

Thus it suffices to prove \lemref{mainlem}
Let $\Phi : \disk \to \Omega$ be a Riemann map,  normalized
so that $S(\Phi)(0) \ne 0$. Conjugate $g$ so that 
$\Phi(0) = \infty$ and $\diam(\partial \Omega) =1$.
If $\Omega/G$ is a surface with punctures, then we can find a
$G$ invariant collection of disjoint  balls $ \cal B_1=\{B^1_j\}$ in $\Omega$, 
each invariant under a parabolic element of $G$ and 
so that $(\Omega \backslash \cup_j B^1_j)/G$ is compact (i.e., 
we are taking a neighborhood of each puncture on $\Omega /G$ and 
lifting it to $\Omega$. Each $B^1_j$ is thus conjugate to one 
of a finite subcollection and each ball has a parabolic 
fixed point of $G$ on its boundary. To each ball in $\cal B_1$
we associate smaller invariant balls $ B^2_j 
\subset B^1_j$ so that the hyperbolic distances between 
$\partial B^1_j$ and $\partial B_2^j $ 
 is  $1$. Also  any  point of 
$\Omega$ which is outside $\cup_j B^2_j$ is within a bounded
hyperbolic distance of the  orbit of $\infty$. Let 
$C_1$ denote this bound. 

If $\Omega/G$ has no punctures, just replace the collections 
$\cal B_j$ by the empty set in the proof that follows. Note that
in this case every point of $\Omega$ is within a bounded distance
of the orbit of $\infty$.

For each $t>0$ consider the  level line of Green's function
in the disk 
$$\Gamma_t = \{ z: |z| = 1-e^{-t} \},$$
and for each orbit point $z_j \in \hat G(0)$ let 
$$\Gamma_t^j = \Gamma_{t(1-|z_j|)} \cap S_j.$$
The first thing we want to see is that the image of 
$\Gamma_t^j$ is very long if $t$ is small enough (independent
of $j$).

\begin{lemma}
For any $M>0$ there is a $t_0$ such that if $t\leq 
t_0$ then  
$$\ell(\Phi(\Gamma_y^j)) \geq M d_j= M \dist(\Phi(z_j), \partial 
\Omega).$$
\end{lemma}

\begin{pf}
This is just \lemref{long-length} after rescaling $z_j$ to 
be the origin.
\end{pf}

We would like to take our orbits in \lemref{mainlem} to 
be (in some sense) taken along the curve $\Gamma_t^j$. 
More precisely, we will break $\Gamma_t^j$ into unit 
hyperbolic segments $\{\gamma_k\}$ and to each segment  associate 
the closest orbit point $z_k$. If there were no parabolic points
then each orbit point would be associated to a bounded
number of segments, say $N$ (depending on only $G$) and  by the 
standard distortion theorems for conformal maps, 
$$ \ell(\Phi(\gamma_k)) \sim d(z_k).$$
Thus by throwing out repeats we would have a collection of 
points $z_k \in S_j$ 
with 
$$ \sum_k d_k \geq \frac MN,$$
and 
$$ d_k \geq C d_j,$$
(where $C$ depends on $t$ and constants in certain distortion
theorems for conformal maps). This proves the \lemref{mainlem}
when there are no parabolics. 

If the surface $\Omega/G$ has punctures, then there may be 
points of $\Gamma_t^j$ which are very far from the closest
orbit of $0$, and we need to replace such pieces by new 
curves which are closer to orbit points.
The idea is that if $\Gamma_t^j$ passes though the ``bottom
half'' of a  horoball 
$B$ in the unit disk, then we can should replace $\Gamma_t^j\cap B$
by the arcs
$$  \Gamma_B= \partial B \cap \{z: \frac 12 t(1-|z_j|) \leq 1- |z|
\leq t(1-|z_j|) \}.$$
These arcs have similar hyperbolic length to the 
the arc $\Gamma_t^j\cap B$, they are only slightly closer
to the boundary of the disk, and they remain within a 
bounded hyperbolic distance of the orbit of $0$.
Most importantly, the Koebe $1/4$-theorem implies that 
$$ \ell( \Phi(\Gamma_B)) \sim \ell(\Phi( \Gamma_t^j\cap B)).$$
Thus we have

\begin{lemma}
For each $z_j$ there is an arc $ \tilde\Gamma_t^j $
(consisting of pieces of $\Gamma_t^j$ and arcs of horoballs such that
\begin{enumerate}
\item  $\tilde \Gamma_t^j \subset S_j \cap \{ |z| \leq 1-\frac t2(1- |z_j|)
\}.$
\item $\ell(\Phi(\tilde \Gamma_t^j)) \geq M d_j$.
\item every component of $\tilde \Gamma_t^j$ has hyperbolic length 
at least $1$.
\end{enumerate}
\end{lemma}

We can now finish the proof of \lemref{mainlem} just as in the 
case without parabolics described above. This finishes the lemma
and hence completes the proof of \thmref{thm1}.

Let us summarize the results from the last three sections.
Suppose $G$ is a finitely generated Kleinian group with a 
simply connected invariant component $\Omega$. 
Recall that $G$ is called an extended Fuchsian group if it is 
Fuchsian or has an index two Fuchsian subgroup.
Then the following are equivalent:
\begin{enumerate}
\item $G$ is not  an extended Fuchsian group.
\item $ \partial \Omega$ is not a circle.
\item $\partial \Omega$ has infinite $1$-dimensional measure.
\item $\dim(\partial \Omega)>1$.
\item $\dim(\Lambda_c(G))>1$.
\item $\delta(G) >1$.
\item $\partial \Omega$ fails to have a tangent somewhere.
\item  $\partial \Omega$ fails to have a  inner tangent and 
all points except possibly  the 
rank 1 parabolic fixed points.
\item The set of inner tangent points of $\partial \Omega$ has 
zero $1$-dimensional measure.
\item Harmonic measure for $\Omega$ is singular to 
$1$-dimensional measure (i.e., there is a subset of $\partial 
\Omega$ of full harmonic measure and $1$-dimensional measure 
zero).
\item Almost every (with respect to harmonic measure) point
of $\partial \Omega$ is a twist point.
\item Harmonic measures for distinct components of $\Omega(G)$
are mutually singular.
\end{enumerate}

 \section{Groups with round components} \label{round} 

Suppose $G$ is  a finitely generated Kleinian group 
and $\Omega(G)$ contains a component $\Omega$ which 
is a  disk.  Then either $\Omega(G)$ consists 
of exactly two components (both round disks) or has 
infinitely many components which
are  disks. In the first case, $G$ is an extended Fuchsian group
 and the limit set
is a circle. The second case is described by the following result.

\begin{thm} \label{circles}
Suppose $G$ is  a finitely generated Kleinian group and 
$\Omega(G)$ contains infinitely many components which
are  disks. Then $\dim(\Lambda(G)) >1$.
\end{thm}

This is slightly different than Larmen's result mentioned
in the introduction, because we do not insist that every 
component be a  disk (it could happen, for example,  that there is 
another component of $\Omega(G)$ with accidental parabolics).
This is known result; it is a special case of Theorem 1 of 
\cite{Canary-Taylor93} which Canary and Taylor prove using a 
result of Furusawa \cite{Furusawa91}. It is also contained in 
results of Sullivan and of Patterson.
Roughly speaking, such a group $G$ must contain 2 Fuchsian 
subgroups $G_1, G_2$ whose limit circles are each contained
in a fundamental domain for the other group. Furusawa's theorem 
then says that $\delta(G_1 * G_2) > \delta(G_1) =1$. Since these
groups are geometrically finite the same inequality holds 
for the dimension of the limit sets.
An alternate approach is to note that the sum of the 
diameters of the orbit of  $\Lambda(G_1)$ under $G_2$  diverges
because $G_2$ is divergence type. Thus
$\Lambda(G_1 * G_2)$ has infinite length. Sullivan's 
results then imply it has dimension $ >1$.

Some very interesting pictures of this type of limit set 
appear in \cite{Bullet-Mantica92}.
Further results on such limit sets and the corresponding 
groups are given in \cite{Keen-Maskit-Series}.

We include a  proof of \thmref{circles} for completeness.
It is a fairly standard computation
involving Hausdorff measures. As usual, we pass
to a finite index subgroup if necessary to remove 
any elliptic elements.
Let $D_1, D_2$ be distinct round components (which are not
just opposite sides of one circle) with stabilizers 
$G_1, G_2$. Since the orbit of $D_1$ under $G_2$ accumulates
densely on $\partial D_2$ and vice versa, we may assume 
(by choosing new disks if necessary) that
$$ \text{dist}(D_1, D_2) \geq 1,$$
$$ \diam(D_1) \leq 1/1000,$$
and the double of each disk in contained in a fundamental 
polygon of the other group (so the translates of the doubles
are disjoint).

Fix values of $\delta >0$ and $N< \infty$. Suppose we 
construct a Cantor set $E$ by an iterative construction
in which a disk $D$ is replaced by at most $N$ disks
$\{ D_j\}$ such that
\begin{enumerate}
\item $D_j \subset 2 D$, and $\{2D_j\}$ are disjoint.
\item $\delta \diam(D) \leq \diam(D_j) \leq \diam(D)/100$.
\item $\sum_j \diam(D_j) \geq 2 \diam(D)$.
\end{enumerate}
Then it is easy to see that $\dim(E) \geq \alpha(\delta, N) >1$.

For $D_1$ we can choose $\delta, N$ and such disks $\{D^1_j\}$
by taking part of the orbit of $D_2$ under $G_1$ (using the 
fact that $G_1$ is divergence type). Similarly for $D_2$ and 
$\{ D^2_j\}$.
At a general step in the construction suppose we have a 
disk $D'$ which is a ``child'' of $D''$ (i.e, $D' \subset
2 D''$ and $\diam(D') \leq \diam(D'')/100$). Then $D'$ corresponds
to either $D_1$ or $D_2$ under the action of  $G_1 * G_2$, so assume it is 
$D_1$. Then there is an element $g \in G_1 * G_2$ so that
$g(D_1) = D'$ and $g(D_2) = D$. Its easy to check that 
$g$ has bounded distortion on $2D_1$ (since it corresponds to the
much smaller disk $D'$). Thus $\{g(D^1_j)\}$ satisfy the 
desired conditions with respect to $D'$ and the construction
may be continued (the constants may be 
different, but we have uniform bounds).

\section {Degenerate limit sets have large $\beta$'s} \label{degen}

Suppose $G$ is a finitely generated degenerate group, 
i.e., $\Omega(G) = \Omega$ has a unique component and 
this component is simply connected. If
 the  Fuchsian equivalent
$\hat G = \Phi \circ G \circ \Phi^{-1}$ has no 
parabolics, then $\Lambda(G)$ has 
``large $\beta$'s'' by the 
argument in \secref{big-Sch}.  We shall show this is true
even if $\hat G$ contains parabolics. The point is that 
a rank 1 parabolic in a degenerate group cannot be doubly 
cusped, and this implies big $\beta$'s.

\begin{thm} \label{deg>1}
If $G$ is a degenerate group there is a $\beta_0=\beta_0(G)>0$
such that if $E = \Lambda(G)$, 
then $\beta_E(Q) > \beta> 0$ for every 
square $Q$ such that $\ell(Q) \leq \text{diam}(Q)$ and 
$ \frac 13 Q \cap \partial \Omega \ne \emptyset$. In particular, 
$\text{dim}(\Lambda(G)) >1$.
\end{thm}


\begin{pf}
 As usual, by Selberg's
lemma we may pass to a finite index subgroup without
torsion (which clearly has the same limit set). Thus we
assume $G$ has no elliptics and let $R = \Omega/ G
= \disk  / \hat G$ 
denote the Riemann surface represented by $G$ and $\hat G$.
If $\hat G$ has no parabolics then  \thmref{deg>1} 
follows from \thmref{Schwarz-dim} as noted in \secref{big-Sch},
so assume $\hat g \in \hat G$ is parabolic and 
let $g = \Phi \circ \hat g \circ \Phi^{-1}$ be the corresponding 
element of $G$. Note that $g$ must also be parabolic 
and call its fixed point $z_0$.

Let $B$ denote a horoball in $\disk$ corresponding to
$\hat g$. Let $\hat \gamma$ denote its boundary and 
$h$ the hyperbolic translation length of $\hat G$ along
$\hat \gamma$. Let $\hat \gamma_0$ denote a segment 
of hyperbolic length $h$ on $\hat \gamma$.
 If we conjugate the group $G$ so that the 
fixed point of $g$ becomes $\infty$ then $g$ becomes a
translation. Then $\hat \gamma_0$ maps to a smooth arc
and thus its easy to see that $\hat \gamma$ maps to a 
smooth quasicircle passing through infinity. In particular, 
its lies in the strip between two parallel lines. 
Thus before conjugating the fixed point to infinity,
we see that $\gamma = \Phi(\hat \gamma)$ is a quasicircle 
passing through $z_0$ and is contained between two 
circles tangent at $z_0$, say $B_1 \subset \Phi(B) 
\subset B_2$. Note that these two balls are invariant 
under $g$ since they correspond to half-planes 
when $z_0$ is conjugated to $\infty$.

By taking  $B_1$ smaller, if necessary, we may assume
that $B_3$, the reflection of $B_1$ through the 
fixed point $z_0$, does not intersect $B_2$.
Then $B_3$ is also invariant under
$g$. We wish to show that $B_3$ must contain points 
of $\Lambda(G)$ that are not too close to $z_0$.
Change variables so that $z_0 =0$, $B_3$ lies in the 
upper half plane and has diameter $1$.  
Let $z = i y$ lie on the positive imaginary axis and note that 
$ |z - g(z)| = O(|y|^2)$, as $y \to 0$.
However, $\hat g$ moves every point of the disk outside
$B$ by at least some hyperbolic distance $\delta$
(depending on $\hat g$ and $B$) and therefore
since $z \in B_3$ lies outside $B_2 \supset \Phi(B)$,
$g$ must move $z$ by more than $\delta$ in the 
hyperbolic metric on $\Omega$. Since $\Omega$ is 
simply connected the hyperbolic metric on $\Omega$
is comparable to $\text{dist}(z,\partial \Omega)^{-1} ds$.
Thus for $z =i y$ small we must have 
$$ \dist (z, \partial \Omega) \leq C \delta^{-1} |z-g(z)|
\leq C \delta^{-1}y^2.$$
By considering the point $z=i y $ with $0< y<\delta/C$ we see that 
$B_3 \cap \Lambda(G)$ must contain a point $z_1$ within
$|y|^2$ of $z$. The orbit of $z_1$ under powers of 
$g$ is a $y^2$-dense set on a circle of radius 
larger that $y/2$ tangent at the fixed point $z_0$.
Since $\Lambda(G)$ is invariant under $g$ it contains 
this sequence of points. 

In fact, we have shown that a point $ z \in B_3 \cap \Omega$
at distance $r$ from $z_0$ can be at most distance 
$ \sim r^2$ from $\partial \Omega = \Lambda(G)$.
In particular, if $Q$ is any square 
so that $Q \cap B_3$ is large (say contains a disk 
or radius $> \ell(Q) /100$) then 
$\beta_{\Lambda(G)} (Q) \geq \beta_0 >0$ (and
$\beta \geq 1/4$ for all sufficiently small square 
containing $0$).

Now we can deduce that $\Lambda(G)$ has large $\beta$'s
everywhere. Suppose we have a
square $Q$ such that $\ell(Q) \leq \text{diam}(Q)$ and 
$ \frac 13 Q \cap \partial \Omega \ne \emptyset$. 
If $\Omega \cap \frac 1{10}Q$ contains no disk of size 
$\ell(Q) /100$, then $\beta_{\Lambda(G)}(Q) \geq 1/20$.
So assume $\Omega \cap Q$ does contain a large disk and
let $z_1 \in \frac 1{10} Q$ be its center. 
If $z_1$ lies outside all the balls $B_1$ corresponding 
to parabolics in $\hat G$ described above, then 
$ w_1 = \Phi^{-1}(z)$ 
lies only a bounded hyperbolic distance from an orbit
where the Schwarzian derivative $|S(\Phi)| (z)(1-|z|^2)^2$
is bounded away from zero. Therefore  
in this case $\beta_{\Lambda (G)} (Q) 
\geq \beta_0$ (where $\beta_0$ depends on $G$, but not 
on $Q$).

Otherwise, $z_1$ must lie inside a horoball $B_1$. Let 
$g$ be the corresponding element of $G$ with fixed 
point $z_0$.  By the choice of $B_1, B_3$ described above
and since $ z_1 \in \frac 1{10} Q$,
 $Q\cap B_3$ must contain a disk of size $\geq \frac 1{100}
\ell(Q)$.  So by our remarks above,
 $\beta_{\Lambda (G)} (Q) \geq \beta_0$ (in fact
it is quite large).
This completes the proof of \thmref{deg>1}.
\end{pf}

One can give yet another proof of $\dim(\Lambda(G)) >1$
for degenerate groups with parabolics using the 
observations above, but without using 
\thmref{big-beta}. One notes that the horoball $B_3$ can 
be packed by an infinite number of ball accumulating to 
the fixed point, each containing points of the limit set 
and so that the sum of the radii diverges. This, together
with the observation that conjugates of $B_1$ take up a 
fixed fraction of the area of $\Omega(G)$ in each such ball,
allows one to define a nested sequence of balls such that the 
resulting Cantor set lies in the limit set and has dimension
larger than one. The packing in question is roughly the 
one corresponding (via a M{\"o}bius transformation)
 to a  unit packing of the upper half plane.
 The dimension of Cantor sets arising from
this kind of packing of the disk were studied in 
\cite{Gardner-Mauldin83}.

\section{Proof of \thmref{thm2} and some corollaries} \label{pf-thm2}

Suppose $G$ has a connected limit set. Then any component
$\Omega$ of $\Omega(G)$ is simply connected. 
The subgroup fixing any component $\Omega$ of $\Omega(G)$ is a
finitely generated Kleinian group  $G_\Omega$ and $\Omega$
is an invariant component of its ordianry set. 
Thus by \thmref{thm1} either $\text{dim}(\Lambda(G))>1$
or every component of $\Omega(G)$ is a   disk. 
If the latter case holds then either $\Omega(G)$ has
two components or infinitely many. If it has two then 
$\Lambda(G)$ is a circle.
Otherwise $\Lambda (G) = \sphere \backslash \cup_j D_j$
for some infinite 
collection of disjoint open disks. 
Thus $\dim(\Lambda(G))>1$ by \thmref{circles}.
(Actually, Larman \cite{Larman67} proved 
that such a set has dimension $>1$, regardless of whether the disks
are associated to any group.)
It follows from the Klein-Maskit combination theorems that
either $\Lambda(G)$ is totally disconnected or $\Lambda(G)$ 
contains a component which is 
itself the limit set of finitely generated subgroup (see
e.g. \cite{Abikoff-Maskit77}, \cite{Maskit88}). Thus
we get 

\begin{thm} \label{9.1}
If $G$ is a finitely generated Kleinian group then its 
limit set is either totally disconnected, a circle or has 
Hausdorff dimension $>1$.
\end{thm}

Our proof of \thmref{thm1} showed that 
if $\Omega$ is a simply connected component which is not a 
disk then $\delta(G) >1$. 
Thus we obtain
\begin{cor} \label{cor-conical}
If $G$ is finitely generated and $\Omega$ is a simply connected 
invariant component of $\Omega(G)$ which is not a 
 disk, then $ \delta(G) = \dim(\Lambda_c(G))>1$.
\end{cor}

It follows from the Klein-Maskit combination theorems that if 
$G$ is finitely generated and $\Lambda(G)$ is totally disconnected 
then $G$ is geometrically finite and hence $\delta(G)=\dim(\Lambda_c(G))= 
\dim(\Lambda(G))$.  Thus  \thmref{9.1} implies 

\begin{cor} \label{dim,delta>1}
If $G$ is a finitely generated Kleinian group then 
$\dim(\Lambda(G)) >1$ iff $\delta(G) >1$.
\end{cor}

Recall that $G$ is called degenerate if $\Omega(G)$ is 
connected and simply connected.
As a special case of \thmref{thm1} we get 
\begin{cor} \label{degenerate}
If $G$ is degenerate then $\dim(\Lambda(G)) >1$
and  $\delta(G) >1$.
\end{cor}

Greenberg \cite{Greenberg66} showed that degenerate groups 
are geometrically infinite so \thmref{gi->dim=2} will imply that
 the limit set of a degenerate
group has dimension $2$, but we mention this weaker 
result here because it is a new result whose 
proof easy,
whereas the stronger result requires much  more.
In \cite{Canary92} Canary shows that $\delta(G)=2$ for 
any degenerate group.

A web group is one in which each component subgroup 
is quasi-Fuchsian (so the ordinary set is a union of 
quasidisks).
By the decomposition theorem of Abikoff and Maskit
\cite {Abikoff-Maskit77} any  finitely generated, 
geometrically infinite group contains either a degenerate
group or a geometrically infinite web group (all component
subgroups quasi-Fuchsian). In either case the 
limit set must have dimension $>1$, so

\begin{cor} \label{gi->delta>1}
If $G$ is a finitely generated, geometrically infinite
Kleinian group then $\dim(\Lambda(G)) >1$ and 
$\delta(G) >1$.
\end{cor}

Later we will show that $\dim(\Lambda(G)) =2$ for 
geometrically infinite groups,  and $\delta(G)=2$
if $\Lambda(G)$ has zero area, but this will be our
best result for  $\delta(G)$ in general.

\def\uMdim{\overline{\text{Mdim}}}
\def\lMdim{\underline{\text{Mdim}}}
\def\Mdim{\text{Mdim}}
\def\Ncov{N_{\text{cov}}}
\def\Cal{\cal}

\def\disk{\Bbb D}
\def\sphere{\overline {\Bbb C}}
\def\dist{\text{dist}}
\def\diam{\text{diam}}
\def\dim{\text{dim}}
\def\Dbar{\overline{\partial}}

\section{Geometry of $\Omega(G)$} \label{G&B-horoballs}

In this section we will prove some lemmas  about the 
geometry of $\Lambda(G)$ and $\Omega(G)$ that we  will use 
later.  The  two main facts we need are
\begin{enumerate}
\item  $\Lambda(G)$ is uniformly perfect.
\item  Show bad horoballs
must be  singly cusped  and far apart.
\end{enumerate}

For $z \in \Omega (G)$ define 
$$ d(z) = \dist( z, \partial \Omega(G)) = \dist (z, \Lambda),$$
where ``distance'' means  spherical distance.
Suppose $\Omega$ is a component of $\Omega(G)$ and let 
$\rho$ denote the hyperbolic metric on $\Omega$. Then a 
standard estimate says that
$$ |d\rho (z) | \leq 2  {|dz| \over d(z)}.$$
Since an orbit is  uniformly separated in the hyperbolic metric this 
implies its is also uniformly  separated in the ``$|dz|/d(z)$''
metric, and hence each Whitney box contains only a bounded
number of orbit points.  Except near parabolic horoballs it is 
also true that Whitney cubes are near at least one orbit 
point. To prove this we need to show $\Lambda(G)$ is 
uniformly perfect.

A compact set $K$ is called {\it uniformly perfect} if there exists
 $\epsilon>0$  such that  for any $ x \in K$ and
$ r< \diam(K)$ there exists $y \in K$ such that
$$ \epsilon r \leq |x-y| \leq r.$$  
There are several well known equivalent formulations of this 
condition (e.g., \cite{Fernandez89}, \cite{Gonzalez92}).
 Suppose $K$ is compact  and $\Omega$ is its 
complement. Then the following are known to be equivalent:
\begin{enumerate}
\item $K$ is uniformly perfect.
\item There is a positive lower bound for the length of the
 shortest closed hyperbolic geodesic in $\Omega$.
\item There is  a constant $C< \infty$ so that 
$$\frac 1C {|dz| \over d(z)} \leq
 |d\rho (z) | \leq 2 {|dz| \over d(z)}.$$
\end{enumerate}

\begin{lemma} \label{lambda=up}
 If $G$ is a finitely generated, non-elementary
 Kleinian group
then $\Lambda$ is uniformly perfect. In particular, if
 $\rho$ is the hyperbolic metric on a component $\Omega$
of $\Omega (G)$, then 
$$ |d\rho (z) | \sim  {|dz| \over d(z)}.$$
\end{lemma}

\begin{pf}
From the discussion above it is enough to verify condition 
(2) for 
$\Omega$: there exists $\epsilon$ so that every closed geodesic
for $\rho$ in $\Omega$ has hyperbolic length $\geq \epsilon$.
By the Ahlfors finiteness theorem, $\Omega(G)/G= R_1 \cup , \dots,
\cup R_N$ is a finite 
union of compact Riemann surfaces with a finite number of punctures, 
so  there is a lower bound  $\epsilon_0$ for the length of the shortest closed 
geodesic (although if there are punctures there is no lower bound
for homotopically non-trivial loops). Suppose there was a 
homotopically non-trivial
closed loop $\gamma$ in $\Omega(G)$ with length less than 
$\epsilon_0$. Then $\Gamma$ projects to a loop in a cusp 
region on $\Omega(G)/G$ which can be homotoped to a puncture.
Thus $\gamma$ can be homotoped in $\Omega(G)$ to curve with 
arbitrarily short hyperbolic length, which implies $\partial 
\Omega(G)$ must have an isolated boundary point. This is 
impossible if $G$ is non-elementary.
Thus  $\Omega(G)$ has
the desired property.
\end{pf}

It is worth noting that the constants in the previous result 
depend on $G$ and cannot be taken to depend, say, only on the 
number of generators. It has also been pointed out to us that 
this lemma (with almost the same proof) was earlier proven 
by Canary in  \cite{Canary91}.

Next we define what we mean by ``good'' and ``bad'' horoballs.
Since $G$ is finitely generated,  the Ahlfors finiteness
theorem \cite{Ahlfors64} says that $\Omega(G)/G$ is a finite 
union of finite  Riemann surfaces $R_1, \dots, R_N$, i.e., each is 
a compact surface with at most a finite number of punctures.
 Let $\{ p_1, \dots, 
p_m\}$ be the punctures in $R = \cup_{i=1}^N r_i$, and for each 
$p_i$ let $B_i^*$ be a neighborhood of $p_i$ which lifts to a 
Euclidean ball $B_i$ in $\Omega$ which is invariant under some 
parabolic element of $G$ (see Lemma 1 of \cite{Ahlfors64}).
Then  $ X=\cup_i R_j \setminus \cup_j B_j^*$ is compact, so  as 
above, we
can choose a finite set of points $E=\{z_1, \dots, z_P\} \subset 
\Omega(G)$ which project to an $\epsilon$  dense  subset of 
$X$.
It will also be convenient to assume (as we may) that each of the 
chosen horoballs $B_j$, is contained in a horoball in $\Omega$
of twice the diameter (so $\partial B_j$ does not come too 
close to $\partial \Omega$ except near the parabolic fixed point.

Since there are only a finite number of surfaces 
in $\Omega(G)/G$, it is easy to see that there is 
an $\epsilon_0 >0$ so that if $\gamma$ is homotopically 
non-trivial loop on one of the surfaces $R_i$ of 
length less than $\epsilon_0$ then $\gamma$ must be 
contained in one of the neighborhoods $B_i^*$. Thus 
any closed curve in $\Omega(G)$ of hyperbolic length 
$\leq \epsilon_0$ is either  homotopically trivial or lies in an 
image of one of the balls $B_i$.

Suppose $B = g(B_i)$ for some
$g\in G$. Suppose $B$ is fixed by a parabolic element
$h \in G$ with fixed point $p \in \partial B$. We say 
that $B$ is doubly cusped it there is a another (disjoint) 
ball $B_1$ fixed by $h$.

Normalize, so that $\Lambda(G)$ has diameter $1$.
  Given a $\eta >0$ we say $B$ is a ``$\eta$-bad''
horoball if
$$ \sup_{z \in \partial B} d(z)
\leq \eta \diam(B).$$
Otherwise we say $B$ is ``$\eta$-good''. We will need the 
following simple facts about  good and bad horoballs.

\begin{lemma} \label{horoballs}
Suppose $\diam(\Lambda)=1$.
\begin{enumerate}
\item There is a $C_1$  (depending only on 
$\eta$) so that for any $\eta$-good horoball $B$, and any $w\in 
\partial B$, there is a point $z \in G(E)$ such that
$$  C_1^{-1}  d(w) \leq 
 d(z) \leq C_1 d(w),$$
$$  C_1^{-1}  d(w) \leq 
|z-w| \leq C_1 d(w).$$
\item There is a $\eta_2$  (depending on $G$) so that if  $B$
is doubly cusped
 then $B$ is $\eta_2$-good.
\item  If $\eta$ is small enough (depending only on $G$), and
if $B$ is  an $\eta$-bad horoball and $D \subset \Omega 
\cap  \eta^{-1/2} B$ is a disk, then $\diam (D) \leq  \eta^{1/6}\diam(B)$.
(the powers are not sharp).
\item For any $\delta>0$ there is a $\eta_3>0$ (depending 
only on $\delta$) such that if $B$ is a $\eta_3$-bad horoball then 
there is  a disk $D \subset 3B$ such that
$ \diam (D) \geq \frac 13 \diam (B),$
and $D \setminus \Lambda$ contains no balls of radius $\geq 
\delta \diam (D)$.
\item There is $\eta_4>0$ so that if  
 $B_1, B_2$ are  $\eta_4$-bad horoballs 
with $\diam(B_1) \leq \diam (B_2)$ then 
$$ \dist (B_1, B_2) \geq 100 \diam (B_1).$$
\item If $B_1, B_2$ are horoballs with
$$\diam (B_1) \leq 
\diam(B_2) \leq 2 \diam (B_1),$$
 and $\dist(B_1, B_2) \leq 
A \diam (B_1),$  then both $B_1$ and $B_2$ are 
$A^{-3}$-good (the power is not sharp).
\end{enumerate}
\end{lemma}

\begin{pf}
The first claim is easy, since if $w \in \partial B$ then 
it is within hyperbolic distance $\epsilon$ of a point in 
$G(E)$, and these points lie in the same or adjacent 
Whitney squares.

To prove (2) suppose $B = g(B_i)$. Then $B$ is doubly 
cusped iff $B_i$ is. Suppose they are and let $\hat B_i$
be the other ball corresponding to $B_i$ and 
let $\hat B = g(\hat B_i)$. Choose $\eta$ so that
$B_i$ is $\eta$-good. If $B$ is $\eta_2$-bad with 
$\eta_2 << \eta$, its easy to see that $g$ maps 
$\hat B_i$ to the exterior of a ball of diameter
$\leq 2 \diam(B)$. But this means that $\diam(\Lambda) 
\leq 2 \diam (B) < 1$, a contradiction. Thus $B$ must be 
$\eta_2$-good, as desired. 

To prove (3), conjugate $B$ by a linear map so it becomes the 
ball of radius $1/2$ centered at $i/2$ and it is 
fixed under an element of the form
$$ h(z) = {z \over 1 + a z}.$$
If $B$ is $\eta$-bad then we must have $a \leq C \eta$.
Therefore, a point $z  \in \eta^{-1/2} B$, is  displaced at most
$$|z-h(z)| = |{a z^2 \over 1+az}| \leq C \eta |z|^2 \leq C \eta^{1/3}
< \eta^{1/6} \epsilon_0/100,$$
if $\eta$ is small enough.
Suppose $z \in  3B \setminus 2B$ is the center of a disk $D \subset
\Omega$ of radius $ \eta^{1/6}$.
Then $D$ intersects
$h(D)$ and connecting the centers gives a circular arc which projects
to a loop on $\Omega(G)/G$ of length less than
$\epsilon_0$. 
Thus iterating $D$ under powers of $h$ gives a circle
fixed by $h$ and lying in $\Omega$.
Moreover, this circle
projects to a curve in $\Omega(G)/G$ of length less than 
$\epsilon_0$. Therefore it lies in a horoball. Hence 
side of this circle must
lie in $\Omega$.  It cannot be the side containing $B$ because this
component contains points of $\Lambda$ (the ones that are within 
distance $\eta$ of $i \in \partial B$). Therefore the other side
must be contained in $\Omega$. This implies $B$ is doubly cusped, 
hence $\eta_2$-good. 
This is a contradiction, so there are no such small disks.

Parts (4), (5), and (6) are all special cases of (3), so we are 
done.
\end{pf}

\def\con{C}

\section{The convex hull of the limit set}\label{convex-hull}

A discrete group of M{\"o}bius transformations is called 
geometrically finite  if there is a finite sided fundamental 
polyhedron for its action on $\Bbb B$. For our purposes, however, 
the following characterization due to Beardon and Maskit
\cite{Beardon-Maskit74} of geometrically finite groups is
much more useful. Recall that a rank 1 parabolic fixed point $p$
is called doubly cusped if there are two disjoint balls in 
$\Omega(G)$, tangent at $p$, and both invariant under the parabolic
subgroup fixing $p$. 

\begin{prop} \label{char.geo.fin}
If $G$ is a Kleinian group then 
$G$ is geometrically finite iff
 $\Lambda(G)$ is the union of $\Lambda_c(G)$, the rank 2 
parabolic fixed points and the  doubly cusped rank 1
parabolic fixed points of $G$.
\end{prop}

Geometrical finiteness can also be characterized in terms 
of the convex hull of the limit set.
If $K$ is a compact set on $S^2 = \partial \Bbb B$ we will 
let $\con(K) \subset \Bbb B$ denote its convex hull with 
respect to the hyperbolic metric on $\Bbb B$. If (as usual in 
this paper) $G$ is a Kleinian group without torsion, 
we let $M= {\Bbb B}/ G$ be the hyperbolic $3$-manifold associated to 
$G$.  Then $C(M) = C(\Lambda(G))/G \subset M$   is called the 
convex core of $M$. Much of the interesting topology of $M$
is associated to the topology of the convex core.

For $r>0$ we define the radius $r$ neighborhood of $C(M)$ as
$$C_r(M) = \{ z \in {\Bbb B}: \dist(z, C(M)) < r \},$$
where distance is  measure in the hyperbolic metric.
Although we don't need it here, we should point out that
Thurston has shown that $G$ is geometrically finite iff 
$C_r(G)$ has finite volume for some (all) $r>0$. 
We cannot take $r=0$  because if $G$ is any Fuchsian 
group, then $\Lambda(G)$ is contained in a circle, so the convex
hull of $\Lambda(G)$ is contained in a hyperplane and hence has 
zero volume. Thus $C(M)$ has finite volume even for infinitely 
generated Fuchsian groups. However, for finitely generated
groups this is not a problem.

\begin{prop} \label{vol->gf}
If $G$ is a finitely generated  Kleinian and $C(M)$ has finite volume
 then $G$ is geometrically finite.
\end{prop}

This is well known, but we will  deduce it as a 
corollary of \lemref{bigC(M)} later
of this section.
The equivalence of the many equivalent formulations
of geometric finiteness is discussed in 
\cite{Bowditch93}.

Our first goal is to see that $C(M)$ can
be separated from $\Omega(G)$ by finite area surfaces. More 
precisely, 

\begin{lemma} \label{surfaces}
Suppose $G$ is finitely generated Kleinian group (without
torsion) and $M$ is the corresponding  hyperbolic $3$-manifold.
Let $\{\Omega_j\}_1^N$ be conjugacy classes of components
of $\Omega(G)$ (i.e., the geometrically finite ends of $M$).
For each $\Omega_j$
There is a  surface $\{S_j\}$ in M  so that 
the following holds. 
\begin{enumerate}
\item $\dist(S_j, C(M)) > 2$.
\item  Each $ S_j $ has finite area.
\item If $U_j= \{ z \in M: \rho(z, S_j) < 1\}$
is a unit hyperbolic neighborhood of $S_j$ then $U_j$ has 
finite hyperbolic volume.
\item The function $f(y) = (\text{vol}(B(y,1))^{-1/2}$ is integrable 
over $U_j$, i.e., 
$$ \int_{U_j} (\text{vol}(B(y,1))^{-1/2} dy < \infty.$$
\item The surfaces $S_j$ separate $C(M)$ from the geometrically 
finite ends of $M$, i.e.,  there is a $\epsilon_0 >0$ so that
$w(z)\leq 1- \epsilon_0$  on the 
component $M_1$  of $M \setminus \cup_j S_j$ containing $C(M)$.
\end{enumerate}
\end{lemma}

{\it Proof:}
This consists of known facts, as described, for example, in 
\cite{Epstein-Marden87}. In fact, we do not even need the 
sophisticated machinery developed there.
 For each component $\Omega_j$ we can 
simply take $S_j$ to be the intersection of a fundamental 
polygon in $\Bbb B$ for $G$ with the surface 
$$ \{z \in {\Bbb B}: 1-|z| = \epsilon_0  \dist(z/|z|, \Lambda(G)\}.$$
By \lemref{lambda=up} (which says $\Lambda(G)$ is uniformly 
perfect)  hyperbolic area on this surface is boundedly 
equivalent to hyperbolic area on $\Omega(G)$ under the radial 
projection if $\epsilon_0$ is small enough (depending only on 
the constants in \lemref{lambda=up}; these depend on $G$, but 
may be taken absolute if $\Lambda(G)$ is connected). Thus (2) 
holds. Condition (3) follows because $U_j$ is boundedly equivalent to 
$S_j \times [0,1]$.  To prove (4) we note that the integral 
obvious converges on any compact piece of $U_j$ and use the
fact that the injectivity radius decreases exponentially to 
bound the integral in the cusps (also note that 
for $y \in U$, $ \text{vol}(B(y,1)) \sim \text{inj}(y)$).
Conditions (1) and (5) are easy to check if $\epsilon_0$ is 
small enough.
$\qed$

We will  also need the Margulis lemma (e.g., \cite{Benedetti-Petronio92}).
 This says that there is an $\epsilon>0$ so that 
if 
$$M_{\text{thin}(\epsilon)} = \{ x \in M: \text{inj}(x) \leq \epsilon \},$$
then every component of $M_{\text{thin}(\epsilon)}$ is 
one of three kinds:
\begin{enumerate}
\item a torus cusp, i.e., a horoball in $\Bbb B$ modulo a 
rank 2 parabolic subgroup.
\item a rank one cusp, i.e., a horoball in $\Bbb B$ modulo a
rank 1 parabolic subgroup.
\item a solid torus neighborhood  of a closed geodesic, i.e., a 
$r$-neighborhood of a geodesic in $\Bbb B$ modulo powers of a
 loxodromic element fixing the geodesic.
\end{enumerate}
Moreover, the components of 
$M_{\text{thin}(\epsilon)}$ (i.e., 
the horoballs  in (1) and (2) and the neighborhoods in (3))
may be taken to be pairwise disjoint. Note that this implies that
if $g \in G$ and $B$ is a ``thin component'' (one of the 
thre types of regions described above), then 
  $g(B)$ hits $B$ iff  $g(B) =B$ and $g$ is in 
the subgroup associated to $B$ above. Actually we will not 
need the full strength of the Margulis lemma. All we will use
is that there is an $\epsilon>0$ as described above, but this 
number may be allowed to depend on $G$.

We will need the following lemma in the proof of 
\thmref{Hdim=lsc}.

\begin{lemma}\label{area>0->big}
If $G$ is a Kleinian group (not necessarily finitely generated)
 and $\Lambda (G)$ has positive
area then $C(M)$ has infinite volume. 
\end{lemma}

Note that we need the hypothesis that $G$ is Kleinian, since
it is possible for $G$ to be a discrete group with $\Lambda (G) 
= S^2$, but $C(M)=M$ to have  finite volume (e.g., if $G$ is 
co-compact).
 
\begin{pf}
We will prove this by showing there is a  sequence $\{ x_n\} 
\in C(M)$ with 
$$\dist(x_n, \partial C(M)) \to \infty,$$
and $\text{inj}(x_n) > \epsilon$ for all n.

For $z \in \Bbb B$ define 
$$ w(z) = \max_{D \subset \Omega(G)} \omega(z, D, {\Bbb B}),$$
where the max is over all round disks in $\Omega(G)$. Then 
$C(\Lambda) = \{ z: w(z) \leq 1/2\}$. Since $w$ is a max of 
harmonic functions it is sub-harmonic but we will not need this.
Instead of $w$ we could simply deal with the distance 
function $dist(z, \partial  C(M))$, which is essentially 
 $ -\log w$.
It is also $G$ invariant, so defines a function on $M$.
It is easy to see that  for any $R>0$ there is an $\epsilon$ so that
$w(x) < \epsilon$ implies $\dist(x, \partial C(M)) > R$.

Since $\Lambda(G)$ has positive area the Lebesgue density 
theorem gives us a point of density $z_0 \in 
\Lambda(G)$. Let $\gamma$ be the hyperbolic geodesic 
connecting the origin to $z_0$ (i.e., a radius of 
$\Bbb B$) and consider points $x \in \gamma$ converging 
to $z_0$. Clearly $w(x) \to 0$, so 
$\dist(x, \partial C(M)) \to \infty$ as $x \to z_0$.
Therefore we only have to show that the injectivity radius 
of $x$ is $\geq \epsilon $ along some subsequence converging to 
$z_0$.

But if this is false then eventually $\gamma$ must be in one 
of the three types of thin regions for all points close enough 
to the boundary. This implies $z_0$ is fixed by some element of 
$G$. Since such points form a countable subset of the limit 
set we may certainly assume $z_0$ is not one of them and we 
are done.
\end{pf}

We will use the following result in the proof of \thmref{gi->dim=2}
(but it is not essential; we will also sketch a proof
which does not require it). We will use part of the proof 
in the proof of \thmref{H=M1} and \thmref{delta=dim1}.

\begin{lemma} \label{bigC(M)}
Suppose $G$ is finitely generated and geometrically infinite.
Then there is a  sequence $\{ x_n\} 
\in C(M)$ with $\dist(x_n, \partial C(M)) \to \infty$ 
and $\text{inj}(x_n) > \epsilon$ for all n.
\end{lemma}

\begin{pf}
The proof is similar to the last lemma but with a few 
extra technicalities because of the possibility  of 
``bad'' horoballs. We will  carry out the proof 
in two  cases
\begin{enumerate}
\item  Every rank 1 parabolic is doubly cusped.
\item There is a rank 1 parabolic which in not doubly cusped.
\end{enumerate}
 Note that case (1) 
contains the case when $G$ has no parabolic elements.

{\it Proof for case (1) of \lemref{bigC(M)}:}
 In this case,
 there is an $\eta>0$ so that every horoball
(if any exist) is $\eta$-good by \lemref{horoballs}.
Let  $z_0 \in \Lambda(G)$ be a point which is 
neither a parabolic fixed point nor in $\Lambda_c(G)$ (such 
a point exists by \propref{char.geo.fin}).
Let $\gamma$ be the hyperbolic
geodesic from the origin in $\Bbb B$ to $z_0$.
We will first show 
$$ \lim_{x\to z_0, x \in \gamma} w(x) =0.$$
Suppose this is false, i.e.,  suppose there is an 
$\epsilon_0$ such that $w(x) > \epsilon_0$ for all 
$ x \in \gamma$. We will prove that  $z_0$ is either in 
$\Lambda_c(G)$ or is a parabolic fixed point.
 
Let  $\{\Omega_j\}_1^N$ be a maximal collection of 
non-conjugate components of $\Omega(G)$. In each we
can find a finite sided fundamental polygon $P_j$ (in $\Omega_j$'s
hyperbolic geometry). Cover $P_j$ by a finite number of 
Euclidean disks, including  a horoball for each parabolic
cusp. Let $ \cal D= \{D_k\}$ denote the (finite) collection of 
disks obtained in this way. Note that there is a constant
$C_0$ (depending only on our choice of $\cal D$) such that
for any disk $D$ in 
$\Omega(G)$  there is a  $D_k$ and a $g \in G$ so that
$ D \cap D_k \ne \emptyset$
 and $\diam(g(D_k)) \geq C \diam 
(D)$. Thus  there is a $C_1 >0$ so that
for any  $x \in \Bbb B$  there is a $D_k$ and $g \in G$ so that 
\begin{eqnarray} \label{omega>w}
 \omega(x, g(D_k), {\Bbb B}) \geq C_1  w(x).
\end{eqnarray}

Therefore we may suppose we have $z_0 \in \Lambda(G)$ 
and  there is an 
$\epsilon_0$ such that for every $ x \in \gamma$ there 
exists $D_k \in \cal D$ and $g\in  G$ such that
$$ \omega(x, g(D_k), {\Bbb B}) \geq \epsilon_0.$$
Let $D(x)= g(D_k)$ be a choice of disk for each $x \in \gamma$.

If $z_0$ is not a parabolic fixed point then $D(x)$ must 
change infinitely often as $x \to z_0$.
Since $D(x)$ is either compactly contained in $\Omega$ or 
a $\eta$-good horoball, 
the proof of \lemref{B3=S2} and \lemref{horoballs} shows that there is
a constant $C_3 >0$ (depending only on our choice of $\cal D$ and 
$\epsilon_0$) so that if 
$$ \omega(x, g(D_k), {\Bbb B}) \geq \epsilon_0,$$
for some $D_k \in \cal C$, then  the orbit point $g(0)$
is contained in the cone with angle $\pi - C_3$ with 
vertex $z_0$. Thus in this case, there is a cone with 
vertex $z_0$ which contains infinitely many images of $0$, 
i.e., $ z_0 \in \Lambda_c(G)$.

Thus we have shown that $w(x) \to 0$ as $x \to z_0$, or 
equivalently $\dist(x_, \partial C(M)) \to \infty$.
The proof now finishes just as in the previous lemma
(if $z_0$ is not a fixed point of some element $g \in G$ then
there is a sequence $\{x_n\} \to z_0$ in $\gamma$ which lies
outside the thin regions). This finishes the proof for case (1).

{\it Proof for case (2) of \lemref{bigC(M)}:}
Now we consider the second case when there is a rank 1 cusp  in 
$\Omega(G)$ which is not doubly cusped. Normalize so that
this parabolic is $z \to z+1$. Since $\infty$ is not doubly 
cusped, we may also assume that $\Omega(G)
\cap \Bbb H$ contains no invariant horoball for this parabolic.
 The corresponding thin part in $\Bbb B$ will be 
$$B_0 = \{(x,y,z): z > C\},$$
 for some $C$.
If $\infty$ is singly cusped, normalize so the lower
half plane is the largest invariant half plane in $\Omega$.
Consider the points $y^* = (0,y,C)$. Since these are 
on the boundary of the thin part the injectivity radii are 
all $\epsilon$. Therefore we need only show $w(y^*) \to 0$
as $y \to \infty$.

Suppose not, say $w(y^*) \geq \nu >0$ and let $D(y)$ be the 
largest disk seen from $y^*$.
Suppose $D(y)$ is not the lower half plane. Then 
 $\diam(D_y)\leq 1$ if 
$y$ is large enough, because otherwise $D(y)$ would not be 
disjoint from its translate under the parabolic element
$z \to z+1$. This would give a line in $\Omega$ and everything
below this line would be in $\Omega$, contradicting that
the lower half plane is the largest half plane in $\Omega$.

If $D(y)$ is constant for all large enough $y$ then clearly 
the harmonic measure of $D(y)$ fromt eh point $y^*$ tends to 
zero, which is what we want.
Otherwise, $D(y)$ must change infinitely often.  
Assume this is the case.

If there is a  sequence $y_n \to \infty$ such that $D(y_n)$ is 
$\eta_n$-bad with $\eta_n \to 0$, then we must 
have $\diam(D(y_n)) \to 0$. This holds because the estimates 
in the proof of \lemref{horoballs} show that for bad horoballs
the associated thin part in $\Bbb B$ has (Euclidean) diameter
which is much greater than the diameter of the horoball in 
$\Omega(G)$. Since  the thin regions in $\Bbb B$ are 
disjoint the thin region corresponding to $D(y_n)$ must 
have diameter $\leq C$ (by the normalization above). Thus
$\diam(D(y_n)) \to 0$. Therefore we are done in this case.

Next suppose that there are no $\eta$-bad horoballs for some 
$\eta>0$. Then we get a subsequence of disks $D(y_n)$ which
are all images of the same disk in $\cal D$ and so that
the diameters converge to some $\nu_0$, $0< \nu \geq \nu_0 \geq 1$.
For two disks in this sequence consider the group element
$g_{nm}$ which maps $D(y_n)$ to $D(y_m)$. Also consider 
a linear conformal map  $h_{n m}$from from $D(y_m) $ to 
$D(y_n)$.  By compactness we
can choose $n,m$ and $h_{nm}$ so that $g_{nm} \circ h_{nm}$
is arbitrarily close to the identity in $\text{PSL}(2, {\Bbb C})$.
Since $h_{nm}$ preserves $B_0$, this contradicts the fact 
(the Margulis lemma) that $g_{nm} (B_0)$ can't intersect $B_0$.
Thus the diameters of $D(y^*)$  must tend to zero.
This completes the proof of \lemref{bigC(M)}.
\end{pf}

Richard Canary has pointed out to us that \lemref{bigC(M)} also
follows easily from known results. His argument goes as 
follows:
 If $C(N)$ is non-compact, then (see \cite{Bonahon86}), there
exists a sequence $\{ \gamma_i\}$ of closed geodesics leaving every compact set.
Pick $\epsilon$ less than the Margulis constant.
Either, $length(\gamma_i) <\epsilon$ or some point on $\gamma_i$ lies in
the $\epsilon$-thick part. In the first case, let $C_i$ denote the component of
the thin part containing $\gamma_i$, then $\partial C_i$ (the boundary
of $C_i)$ intersects $C(N)$ (otherwise $C(N)\subset\partial C_i$
which implies that $N$ is elementary), so pick $x_i$ in the
the intersection of $\partial C_i$ and $C(N)$.
In the second case pick $x_i\in \gamma_i\cap N_{thick}$. Then,
$\{ x_i\}$ leaves every compact set and every point in the
sequence has injectivity radius at least $\epsilon$.

Note that we obtain \propref{vol->gf} as a corollary of 
\lemref{bigC(M)} since the lemma gives us infinitely 
many disjoint balls which all have volume bounded uniformly 
from below. The proof also shows the following known fact (e.g., 
\cite{Maskit88}, it also follows from \propref{char.geo.fin}):

\begin{cor}\label{gf->dc}
If $G$ is a finitely generated, geometrically finite group, 
then every parabolic cusp of $\Omega(G)$ is doubly cusped
(hence $\eta$-good for some $\eta$ depending on $G$).
\end{cor}

We will use this in the proof of \thmref{H=M1} and \thmref{delta=dim1}.

\newpage 

\section{Geometrically infinite groups} \label{geo-inf}

In this section we will prove, 

\begin{thm} \label{dim=2}
If $G$ is a finitely generated and geometrically infinite 
then 
$$\dim(\Lambda(G)) =2.$$
\end{thm}
 
If $\delta(G) =2$ then this follows from \thmref{delta=conical1}.
Therefore we may assume $\delta(G)<2$.  In this case
\thmref{dim=2}  follows from 

\begin{thm} \label{area>0}
If $G$ is a finitely generated, geometrically infinite 
group and $\delta(G) < 2$ then $\Lambda(G)$ has positive
area.
\end{thm}

Note that the Ahlfors conjecture and \thmref{area>0} imply that
$\delta(G)=2$ for any finitely generated, geometrically infinite
group $G$.

{\it Proof of \thmref{area>0}:}
Let $\Lambda= \Lambda(G)$ be the limit set of $G$. As usual, we assume 
$G$ has no torsion. Let $\delta = \delta(G)$ be the 
critical index for the Poincar{\'e} series and 
$\lambda_0$ the base eigenvalue for the Laplacian 
on $M = {\Bbb B} / G$. By \corref{gi->delta>1} we know that
if $G$ is geometrically infinite then $\delta(G) >1$.
Therefore by a result of Sullivan (Theorem 2.18
\cite{Sullivan87})
$$ \lambda_0 = \delta(G)(2-\delta(G)).$$

Let $C(\Lambda)$ be the convex hull of the limit set in 
$\Bbb B$ and $C(M) = C(\Lambda)/G \subset M$ be the convex 
core of $M$ (see previous section). 
Let $S = \cup_j S_j$  and  $U = \cup_j U_j$ be the surfaces and 
their unit neighborhoods given by \lemref{surfaces}.
Let $M_1$ be the component of $M \setminus \cup_j S_j$.
We will prove that Brownian motion is transient in $M_1$. 
More precisely, we will show that

\begin{lemma} \label{transience}
Suppose $M, S, M_1$ are as above.
Then there is a point 
$x_0 \in M_1$ such that a Brownian motion started at 
$x_0$ has a positive probability of never hitting  $S$.
\end{lemma}

First let us see why this implies $\Lambda$ has positive area.
Let $\tilde x$ be a point in the hyperbolic ball $\Bbb B$ 
which corresponds to the point $x_0$ in the previous paragraph
and let $\tilde M_1$ be the lift of $M_1$ to the ball. 
Then a Brownian motion started at $\tilde x_0$ has a positive
probability of hitting $S^2 = \partial \Bbb B$ without leaving
$\tilde M_1$. Since $w< 1- \epsilon_0$, such a path must 
first hit $S^2$ at a point of $\Lambda$. Since harmonic 
measure for $\Bbb B$ is mutually absolutely continuous 
with respect to area measure on $S^2$ this means that 
$\Lambda$ has positive area.
Therefore we need only prove \lemref{transience}.

{\it Proof of \lemref{transience}:}
The lemma follows from known estimates on the heat kernel
$K(x,y,t)$ for a Riemannian manifold.
 

\begin{prop}[Davies, Theorem 17 \cite{Davies88}] \label{davies}
If $\delta >0$ then there is a constant $C = C(\delta)$ such that
$$0 \leq K(x,y,t)  \leq C \text{vol}(B(x,1))^{-1/2} \text{vol}
(B(y,1))^{-1/2} e^{(\delta- \lambda_0)t} e^{-\rho(x,y)^2
(4+\delta)^{-1}t^{-1}},$$
where $\rho(x,y)$ denotes the hyperbolic distance between $x$ and
$y$ in $M$.
\end{prop}

See also \cite{Davies93} and \cite{Grigor'yan93} and their references.
To apply Davies' result we let $\delta = \lambda_0/2$.
Choose $R>0$ very large.
By \lemref{bigC(M)}  we can choose a point $x \in C(M)$ 
where the injectivity radius $\text{inj}(x) > \epsilon >0$ 
is uniformly bounded below and 
$$ \dist(x, U)\geq \dist(x, \partial C(M)) \geq R.$$
The expected time a Brownian motion started at $x$ spends in $U$ is 
$$  \int_0^\infty [\int_U K(x,y,t) dy] dt .$$
By  Davies' estimate this is bounded by 
$$ \text{vol}(B(x,1))^{-1/2}  \int_0^\infty
 e^{- \lambda_0t/2} e^{-R^2/5t} dt \cdot 
 \int_U \text{vol}(B(y,1))^{-1/2}dy .$$
Using the fact that $\text{inj}(x)$ is bounded below and 
 part (4) of \lemref{surfaces} this is bounded by 
$$C \int_0^\infty 
 e^{- \lambda_0t/2} e^{-R^2/5t} dt.$$
This can be estimated by breaking the integral at $t=R$ to get 
$$C \int_0^R e^{- \lambda_0t/2} e^{-R/5} dt
+C \int_R^\infty e^{- \lambda_0t/2}  dt
\leq  \frac C{\lambda_0} (e^{-\frac 15 R} + e^{- \frac {\lambda_0}2 R}).$$
Thus by choosing $R$ sufficiently large we can make the 
 the expected time a Brownian motion started
at $x$ spends in $U$ is a small as we wish. 

From this we wish to  deduce that a Brownian motion started at $x_0$ 
never hits $S$ with some positive probability. To do this we simply 
note that  the expected time it takes a Brownian motion started
 at a point $y$  of $S$  to first leave $U$  (i.e., to travel unit 
distance from $S$) is bounded away from zero independent of the 
starting point $y$. This is because the expected time to travel 
distance $1$ in $M$ is greater than or equal the expected time to 
travel this distance in the covering space $\Bbb B$, and this is 
bounded away from zero. Let $t_0>0$ be a lower bound for the 
expected time to travel distance $1$.

Thus the expected time a Brownian motion started at $x$ spends in 
$U$ is at least the probability that it every hits $S$ times 
the bound $t_0$. Now choose $x_1$ so that
the expected time spent in $U$ is less than $t_0/2$. Then the 
probability that a Brownian motion starting at $x_1$ ever hits
$S$ is $< 1/2$. This completes the proof of \thmref{dim=2}.$\qed$

We do not actually need the full strength of \lemref{bigC(M)}.
All we need is that there is a sequence $\{x_n\}$ with
$\text{inj}(x_n) \geq \epsilon$ and $x_n$ leaving every 
compact set (rather than $\dist(x_n, \partial C(M))\to \infty$).
This is because we can cut $U$ up into a compact piece and 
a finite number of small tips of cusps. If the volume of the 
small cusp pieces are small enough, then any uniform bound
on the heat kernel shows that the expected time in these 
pieces is small. The time in the remaining compact part 
is small as long as the distance of $x_n$ form the compact
part is large. This is true as long as $x_n$ tends to 
infinity in $M$.

{\it Sketch of second proof of \lemref{transience}:}
It actually suffices to use a weaker estimate on the heat
kernel and less information on the convex hull. All we need
to know is that Brownian motion is transient on a complete, connected, 
infinite volume Riemannian manifold with lowest eigenvalue
bounded away from zero.   For example, the estimate
$$ K(x,y,t) \leq C_1 e^{-C_2 t},$$
for some fixed $x$, $y$ in some ball and $t> T_0$ would 
be sufficient for this.

We pay for the less precise estimate by a more involved 
construction on the manifold. Let $S$ be the surfaces
described by \lemref{surfaces} which separate the convex core
from the geometrically finite ends. Cut $M$ along $S$ and 
let $M_1$ be the component containing $C(M)$. Glue
two copies of $M_1$ along $S$. We claim that the resulting 
manifold $N$ (the double of $M_1$)
 has lowest eigenvalue bounded away from zero.
If so then the heat kernel estimates apply to the new manifold
$N$ and we deduce that  the expected time a  Brownian motion 
spends in $U$ (the unit radius neighborhood of the $S$) is 
finite. By the Borel-Cantelli lemma this says that there are 
points $x$ in $M$ from which the probability of ever hitting 
$S$ is strictly less than $1$ (in fact is as small as we wish).
Thus there is a point $x \in N \setminus U$ from which the 
chance of ever hitting $U$ is less than $1/2$. But the 
two components of $N \setminus U$ are both exactly $M_1$. 
Thus Brownian motion in $M_1$ must have a positive probability of 
tending to infinity without ever hitting $S$.

 This proves the lemma, except for verifying that $N$ has first
eigenvalue bounded away from $0$. We will not verify this in detail, 
but simply note that since $M$ has constant negative curvature and 
lowest eigenvalue $>0$, Buser's inequality (e.g. \cite{Canary92})
implies the Cheeger constant for $M$ is bounded away from zero. 
From this one proves that the Cheeger constant for the manifold
with boundary $M_1$ is non-zero, and from this that the Cheeger 
constant for the doubled manifold $N$ is non-zero. Then Cheeger's
estimate say the first eigenvalue  for $N$ is non-zero, as desired.
$\qed$

\section{Lower semi-continuity of Hausdorff dimension}\label{LSC-dim}

In this section we will prove, 

\begin{thm} \label{dim=lsc}
If $G$ is a finitely generated M{\"o}bius  group and $\{G_n\}$
is a sequence of  M{\"o}bius groups converging algebraically to $G$
then 
$$ \liminf_n \dim(\Lambda(G_n)) \geq \dim(\Lambda(G)).$$
\end{thm}

In particular, if $G$ is geometrically infinite then this result 
and \thmref{gi->dim=2} imply that $\lim_n \dim(\Lambda(G_n))=2$.

If $G$ is not discrete then \thmref{dim=lsc}
 follows from \thmref{delta=conical1} and the remarks
following \corref{delta=lsc}. If all but finitely many 
of the $\{G_n\}$  are non-discrete then by passing to a 
subsequence we may assume $\Lambda(G_n)$  is always all 
of $S^2$,  are all circles or are always two or less points.
We may also assume that the limit sets $\Lambda_n=\Lambda(G_n)$ converge
in the Hausdorff metric to a set $\Lambda_\infty$. Moreover, 
because of the special form of the sets we must have 
$$ \dim(\Lambda_\infty) \leq  \liminf_n \dim(\Lambda_n).$$
It is easy to check that the limit set of $G$
is a subset of the limit of the 
sets $\Lambda_n$.

Therefore we may assume all the groups are discrete.
If all but finitely many of the $\{G_n\}$ satisfy 
$\Lambda(G_n) = S^2$, there is nothing to do, so we may as
well assume (by passing to a subsequence) that the 
$\{G_n\}$ are all Kleinian groups.

If  $\Lambda (G)$ has zero area or $\delta(G) =2$
 then $\delta(G) = \dim(\Lambda(G))$
and this follows from \corref{delta=lsc}. 
 Therefore we need only 
consider the case when $\delta(G) < 2$ and 
 $\Lambda(G)$ has positive area (possibly the whole 
sphere)
If $G$ where geometrically finite then either $\Lambda(G)=0$
(\cite{Ahlfors66}) or $\delta(G)=2$.  In either case we are finished, 
so we may assume that $G$ is geometrically finite.

By passing to a subsequence we may assume 
$$  \lim_n \dim(\Lambda(G_n)),$$
exists and equals the liminf of the original sequence.
If all but finitely many of the groups $\{G_n\}$ are geometrically 
infinite then
$$ \liminf_n \dim(\Lambda(G_n)) =2,$$
by \thmref{dim=2} so there is nothing to do. Otherwise, by 
passing to another subsequence, we may assume
all the $G_n$'s are geometrically finite.
We may also assume
$ \lim_n \delta(G_n)= \alpha  $ exists and is strictly less 
than 2 (otherwise the result follows from the inequality 
$\dim(\Lambda(G_n)) \geq \delta(G_n)$).
To finish the proof we need

\begin{thm} \label{delta->2}
If $\{G_n\}$ is a sequence of geometrically finite Kleinian  groups which 
 converges algebraically to a finitely generated, 
geometrically infinite  discrete group $G$ then $\delta(G_n ) \to 2$.
\end{thm}

This result  follows from two known results

\begin{prop}[Canary] \label{lambda<vol}
If $G$ is a $n$-generated,  geometrically finite group  then 
$$ \lambda_0 \leq  { A_n  \over \text{vol}(C(M)) },$$
where $A_n$ is a constant that only depends on the number of 
generators of $G$.
\end{prop}

{\it Proof:}
This  is essentially Theorem A of \cite{Canary92} except that 
there Canary proves 
$$ \lambda_0 \leq A  {   \chi(\partial C(M)) \over \text{vol}(C(M)) },$$
where $A$ is an absolute constant and 
where $\chi$ denotes the Euler characteristic. However, the 
Euler characteristic of $\partial C(M)$ is the same as that 
of $\Omega(G)/G$, because there is always a homeomorphism 
between the two (e.g., see Epstein and Marden's paper 
\cite{Epstein-Marden87}).
By the Bers inequality (a quantitative version of the Ahlfors 
finiteness theorem) the area, and hence the Euler characteristic,
 of $\Omega(G)/G$ can be bounded in terms of $n$, the number of 
generators of the group $G$.
Thus Canary's results says 
$$ \lambda_0 \leq  { A_n   \over \text{vol}(C(M)) },$$
where $A_n$ depends only on the number of generators.
$\qed$


\begin{prop} \label{boundedvol->gf}
Suppose $\{G_n\}$ is a sequence of $n$-generated, geometrically 
finite  Kleinian groups such that 
$$ \sup_n \text{vol}(C(M)) \leq M < \infty.$$
If the sequence $\{G_n\}$ converges algebraically to a finitely generated, 
geometrically infinite,
discrete group $G$  then $\Lambda(G)$ has zero area.
\end{prop}

{\it Proof:}
This is an easy case of a result obtained by E. Taylor in 
\cite{Taylor94} and is probably well known.
 Here we will sketch a proof which follows 
an argument given by J{\o}rgensen and Marden in 
\cite{Jorgensen-Marden90}.

Suppose $\Lambda(G)$ has positive area. We will derive a 
contradiction.
 Suppose $\{G_n\}$ is a 
sequence of $N$-generated Kleinian groups converging
algebraically to a Kleinian group $G$.
 The Borel selection
theorem says that the set of compact subsets of a compact 
metric space is itself compact with the Hausdorff metric
$$d(E,F) = \max_{z\in E} \dist(z, F) + \max_{w\in F} \dist(w, E),$$
so by passing to a subsequence (which we also denote $\{G_n\}$)
we may assume the sets $\Lambda_n = \Lambda(G_n)$ converge in the
Hausdorff metric to a compact set $\Lambda_\infty$. 

We say a sequence  $\{G_n\}$ converges polyhedrally to a 
group $H$ if   $H$ is discrete and for some $x_0 \in \Bbb B$,
the fundamental polyhedra (the Dirichlet polyhedron)
$$P(G_n)= \{ z \in {\Bbb B}: \rho(z,x_0) \leq \rho(z, g(x_0)) \text{ for all }
g \in G_n \},$$
converge to $P(H)$ uniformly on compact subsets of $\Bbb B$.
By Proposition 3.8 of \cite{Jorgensen-Marden90} any algebraically 
converging subsequence has a polyhedrally convergent  subsequence
and the polyhedral limit contains the algebraic limit (but they 
need not be equal).

A third notion of convergence of groups is geometric convergence.
Given a sequence of groups $\{G_n\}$  we define 
$$ \text{Env}\{G_n\} = \{ g \in \text{PSL}(2, {\Bbb C}): g 
= \lim _n g_n,  g_n \in G_n \},$$
and we say $G_n \to H = \text{Env}\{G_n\}$ geometrically if 
for every subsequence $\{ G_{n_j}  \}$, $\text{Env}\{G_{n_j}\}
=\text{Env}\{G_{n}\}$. 
Proposition 3.10 of \cite{Jorgensen-Marden90} says that $G_n$ converges
geometrically to $H$ if and only if it converges polyhedrally 
to $H$.

Thus we may assume that we have groups $G \subset H$ such that 
\begin{enumerate}
\item $ G_n \to G$ algebraically,
\item $ G_n \to H$ polyhedrally and geometrically,
\item $ \Lambda_n \to \Lambda_\infty$ in the Hausdorff metric.
\end{enumerate}

We first claim that $\Lambda(H) \subset \Lambda_\infty$.
If $\Lambda_ \infty = S^2$ there is nothing to do, so we may 
assume $\Lambda_\infty$ is not the whole sphere.
In this case we follow the proof of Proposition 4.2 of 
\cite{Jorgensen-Marden90}. Let $\Omega = S^2 \setminus \Lambda_\infty$
and suppose $K, K'$ are compact sets such that
$$ K \subset \text{int}(K') \subset K' \subset \Omega.$$
Suppose $h \in H$.  Because $\{ G_n\}$ converges 
geometrically to $H$ we can write $h =\lim_n g_n$ with 
$g_n \in G_n$.

We claim that $h(K) \subset \Omega$. If not then $h(K)$
intersects $\Lambda_\infty$, so $\text{int}(h(K')$ also 
hits $\Lambda_\infty$. This implies that $\text{int}(h(K')$
intersects $\Lambda_n$ for all large enough $n$, say $n \geq N_1$.
Therefore $\text{int}(g_n(K')) $ hits  $\Lambda_{m}$ for all 
$m \geq N_1$, for all sufficiently large $n$, say $n\geq N_2$. So if 
$N_3 = \max (N_1, N_2)$, then $n \geq N_3$ implies 
$$  \text{int}(g_n(K')) \cap \Lambda_n \ne \emptyset.$$
Therefore
$$  \text{int}(K') \cap  g_n^{-1}(\Lambda_n) =
\text{int}(K')\cap  \Lambda_n \ne \emptyset.$$
This is a contradiction, so we must have $h(K) \subset \Omega$.
This implies $h(\Omega) \subset \Omega$. Since the same 
argument applies to $h^{-1}$, we see that $h(\Omega) = \Omega$,
or equivalently, $h(\Lambda_\infty) = \Lambda_\infty$.
Since $\Lambda_\infty$ is a closed set which is invariant under 
the group $H$ we must have $\Lambda(H) \subset \Lambda_\infty$ as 
desired (recall that the limit set is the smallest closed $H$-invariant
set if $H$ is non-elementary).

Since $\Lambda(H) \subset \Lambda_\infty$ the convex hull 
$C(\Lambda(H))$ of $\Lambda(H)$ in $\Bbb B$ is contained in the convex 
hull  $C(\Lambda_\infty)$ of  $\Lambda_\infty$. The convex hulls
$C(\Lambda_n)$ of the sets 
$\Lambda_n$ converge, uniformly on compacta, to 
$C(\Lambda_\infty)$. Thus for any  $R< \infty$,
$$\text{vol}(C(\Lambda(H)) \cap P(H) \cap B(x_0, R))
\leq \liminf_n  \text{vol}
(C(\Lambda_n) \cap P(G_n) \cap B(x_0, R)) 
\leq M.$$
Thus $C(\Lambda(H)) \cap P(H)$ has finite volume.
Here we are using the convexity of the sets $\Lambda_n$; in 
general, sets of zero volume can converge in the Hausdorff metric
to sets of infinite volume.

The volume of the quotient manifold is easily seen to be 
continuous under geometric convergence, so we deduce
that $ {\Bbb B} / H$ has infinite volume.
If $ \Lambda(H)= S^2$ then $\text{vol}(M) = \text{vol}(C(M)$ is finite, 
which is a contradiction. Thus 
$\Lambda(H) \ne S^2$.


If $\Lambda(H)$ is a proper subset of $S^2$ it still has positive
area since it contains the set $\Lambda(G)$.
Thus by \lemref{area>0->big}
the convex core of $H$ has infinite volume. This is another contradiction 
and so completes the proof of \propref{boundedvol->gf}.
$\qed$

Now that we have the two propositions, we can 
can finish the proof of \thmref{delta->2}. Suppose 
$G$ is a finitely generated, geometrically infinite discrete group and 
$\{G_n\} $ are geometrically finite groups converging to 
$G$  algebraically.
If $\Lambda(G)$ has zero area then $\delta(G) = \dim(\Lambda(G))=2$ and 
$\delta(G_n) \to 2$ by \corref{delta=lsc}.
Thus we may assume $\Lambda(G)$ has positive area.
  By \propref{boundedvol->gf}
we must have 
$$ \text{vol}(C(G_n)) \to \infty,$$
so by \propref{lambda<vol} we get 
$$ \lambda_0(G_n) \to 0.$$
Since $G$ is geometrically infinite we have 
$\delta(G) >1$ by \corref{gi->delta>1} and so by 
\corref{delta=lsc} we have 
$$\liminf_n \delta(G_n) \geq \delta(G) > 1.$$
Therefore we may as well assume $\delta(G_n)\geq 1$
for all $n$.
By Sullivan's result \cite{Sullivan87} $\lambda_0 = \delta(2-\delta )$  if
$\delta\geq 1$ 
so,
$$ \delta(G_n) \to 2,$$
as $n \to \infty$. This completes the proof of 
\thmref{delta->2}. $\qed$

\section{Upper Minkowski dimension} \label{UMdim=BT}

In this section we introduce upper Minkowski dimension and 
prove a simple lemma which we will need in the next section.

Suppose $K$ is a bounded  set in ${\Bbb R}^d$ (or any metric
space for that matter) and let $N(K, \epsilon)$ be the minimal
number of $\epsilon $ balls needed to cover $K$. We define the
{\it upper Minkowski dimension} as
$$ \uMdim (K) = \limsup_{\epsilon \to 0} { \log N(K , \epsilon )
\over \log 1/\epsilon},$$
and the {\it lower Minkowski dimension}
$$ \lMdim (K) = \liminf_{\epsilon \to 0} { \log N(K , \epsilon )
\over \log 1/\epsilon}.$$
If the two values agree, the common value is simply called the
Minkowski dimension of $K$ and is denoted $\Mdim(K)$.
 The Minkowski dimension says the
the number of balls of size $\epsilon $ needed to cover
$K$ grows  like  $\epsilon^{-\Mdim(K)}$.
It is clear that we always have  
$$ \dim(K) \leq \lMdim(K) \leq \uMdim(K),$$
and strict inequality is possible in both places.
   A simple example where
the Minkowski dimension exists, but is strictly larger than
the Hausdorff dimension, is given by the countable set 
$ K = \{0\} \cup \{ \frac 12, \frac 13, \frac 14, \dots\}$.
Here it is not to difficult to prove
$$ N(K, \epsilon) \sim \frac 1{\sqrt \epsilon},$$
so that $\Mdim(K) = 1/2$. However, since the set is countable,
$\dim(K) =0$. A much more sophisticated class of sets
where the two dimensions disagree is given in  \cite{McMullen84}.
 Given a compact set $K$ in the plane let $\Omega =
{\Bbb R}^2 \setminus K$ be its complement. A Whitney decomposition of
$\Omega$ is a collection of squares $\{ Q_j\}$ which are disjoint,
except along their boundaries, and such that
$$ \frac 1{10} \dist(Q_j, \partial \Omega) \leq \ell(Q_j)  \leq
10 \dist(Q_j, \partial \Omega).$$
The existence of Whitney decomposition for any open set is
a standard fact in real analysis (e.g., \cite{Stein71}).
The squares $\{Q_j\}$ may also be taken to be dyadic.
 
For any compact set $K$ we can  define an exponent  of convergence
similar to the exponent of convergence of a Poincar{\'e} series,
$$ \kappa = \kappa(K) = \inf \{ \alpha: \sum_j \ell(Q_j)^\alpha
< \infty\},$$
where the sum is taken over all squares in Whitney decomposition of
$\Omega = K^c$ which are within distance $1$ of $K$ (we have to
drop the ``far away'' squares or the series might not converge).
It is easy to check that this does not depend on the particular
choice of Whitney decomposition.
This number has been rediscovered many times in the literature,
but seem to have been first used by Besicovitch and Taylor
in \cite{Besicovitch-Taylor}. See also  Tricot's paper\cite{Tricot81}
where he  shows
that it agrees with the upper Minkowski dimension. We include
a proof of this fact since it is central to our argument.
 
\begin{lemma} \label{kappa=uMdim}
For any compact set $K$, $\kappa\leq \uMdim(K)$. If $K$ also has
zero area then $\kappa = \uMdim(K)$.
\end{lemma}
\begin{pf}
We start with the easy direction, $\kappa \leq \uMdim(K)$.
Let $D = \uMdim(K)$ and choose  $\epsilon >0$. Then by the definition of
$\uMdim(K)$ we have that 
$$ N(K, 2^{-n}) \leq C 2^{D+\epsilon}.$$
Let $\Cal S_n$ be a covering of $K$ by  fewer than 
$C 2^{D+\epsilon}$ squares of size $2^{-n}$. If $Q$ is a
Whitney square  with $ 2^{-n-1} \leq \ell(Q) < 2^{-n}$, then
choose a point $x \in K$ with $\dist(x, Q) \leq \ell(Q)$. Let
$S(Q) \in \Cal S_n$ be the square containing $Q$. Since $S(Q)$ and
$Q$ have comparable sizes and there distance apart is at most
$\ell(Q)$, we easily see that each $S \in \Cal S_n$ can only be
associated to a uniformly bounded number of $Q$'s in the Whitney
decomposition, say $A$.
 
Let $W_n$ be the number of Whitney squares with size 
 $ 2^{-n-1} \leq \ell(Q) < 2^{-n}$.
Then 
\begin{eqnarray*}
  \sum_j \ell(Q_j)^{D + 2 \epsilon} 
&\leq& \sum_{n=0}^\infty  W_n 2^{-n(D+2 \epsilon)} \\
&\leq& \sum_{n=0}^\infty  C N(K,2^{-n}) 2^{-n(D+2 \epsilon)} \\
&\leq& C \sum_{n=0}^\infty  2^{n(D+\epsilon)} 2^{-n(D+2 \epsilon)} \\
&\leq& C \sum_{n=0}^\infty  2^{-n\epsilon} \\
&<& \infty,
\end{eqnarray*}
which proves $\kappa \leq D + 2 \epsilon$. Taking $\epsilon \to 0$
gives $\kappa \leq D = \uMdim(K)$.

Now we assume $K$ has zero area and will prove 
$\kappa \geq D= \uMdim(K)$. As above, let $\epsilon >0$ and
suppose $\{ Q_j\}$ is a Whitney decomposition of $ \Omega
= K^c$. By the definition of $\uMdim(K)$ we have 
$$ N(K, 2^{-n}) \geq 2^{n(D- \epsilon)},$$
for infinitely many $n$. Suppose $n_0$ is a value where this occurs
and let $\Cal S = \{ S_k\}$ be a covering of $K$ with dyadic squares
of size $2^{-n_0}$. Let for each $S_k \in \cal S$,
 $\Cal C_k= \{ Q_{j k}\}$ be the collection of 
Whitney squares which intersect $S_k$. If we assume the $Q_j$ are
dyadic, then every square hitting $S_k$ is contained in $S_k$.
 Since the area of $K$ is zero,  this gives 
\begin{eqnarray*}
2^{-2n_0} =\text{area}(S_k) = \text{area}(S_k \setminus K)
= \text{area}(S_k \cap \Omega) = \sum_{\Cal C_k}
\ell(Q_{j k})^2.
\end{eqnarray*}
Therefore, 
\begin{eqnarray*}
\sum_{\Cal C_k} \ell(Q_{j k})^{D-2\epsilon} 
&=& \sum_{\Cal C_k} \ell(Q_{j k})^2 \ell(Q_{j k})^{-2+D - 2\epsilon}\\
&\geq& \sum_{\Cal C_k} \ell(Q_{j k})^2 \ell(S_k)^{-2+D - 2\epsilon}\\
&=&  \ell(S_k)^{-2+D - 2 \epsilon}\sum_{\Cal C_k} \ell(Q_{j k})^2 \\
&=&  \ell(S_k)^{-2+D - 2 \epsilon} \ell(S_k)^2 \\
&=&  \ell(S_k)^{D - 2 \epsilon} \\
&=&  2^{-n_0(D - 2 \epsilon)}. \\
\end{eqnarray*}
Hence,   
\begin{eqnarray*}
\sum_{j} \ell(Q_j)^{D-2 \epsilon}
&\geq &  \sum_{k} \sum_{\Cal C_k} \ell(Q_{j k})^{D - 2 \epsilon}\\
&\geq &  \sum_{k}\ell(S_k)^{D - 2 \epsilon}\\
&\geq &  N(K, 2^{-n_0})2^{-n_0(D - 2 \epsilon)}\\
&\geq &  2^{n_0(D - \epsilon)} 2^{-n_0(D - 2 \epsilon)}\\
&= &  2^{n_0\epsilon} .
\end{eqnarray*}
Taking $n_0 \to \infty$, we get
$\sum_{j} \ell(Q_j)^{D-2 \epsilon} = \infty$, and hence
$\kappa \geq D- 2\epsilon$. Taking $\epsilon \to 0$ gives the
desired result.
\end{pf}

\section{ Minkowski dimension equals Hausdorff dimension}\label{Minkowski}

In this section we prove, 

\begin{thm} \label{H=M}
If $G$ is a finitely generated Kleinian group then the 
Minkowski dimension of $\Lambda$ exists and equals the 
Hausdorff dimension.
\end{thm}

Finitely generated is necessary, because it is possible to 
construct infinitely generated Kleinian groups where the 
Minkowski dimension of $\Lambda(G)$ fails to exist.

\begin{thm} \label{delta=dim}
If $G$ is finitely generated Kleinian group and $\text{area}(\Lambda(G))
=0 $ then $\delta(G) = \dim(\Lambda)$.
\end{thm}

Thus if the Ahlfors conjecture holds, we would have 
$\delta(G) = \dim(\Lambda(G))$ for all finitely generated
groups.
\thmref{H=M} and \thmref{delta=dim} follow from this
 more technical looking 
result.

\begin{thm} \label{delta=uM}
Suppose $G$ is a finitely generated Kleinian group.
If $\text{area}(\Lambda(G))=0$ then $\delta(G) = 
\uMdim(\Lambda (G))$.
\end{thm}

{\it Proof of \thmref{H=M} and \thmref{delta=dim}:}
To see how \thmref{delta=uM} implies \thmref{H=M} we consider 
two cases. First, if $\dim(\Lambda)=2$ then 
$$2= \dim(\Lambda) \leq \lMdim(\Lambda) \leq \uMdim(\Lambda) \leq 2,$$
so all are equal to $2$. On the other hand, if $\dim(\Lambda) < 2$,
then $\Lambda$ has zero area, so \thmref{delta=uM} applies and 
gives
$$\dim(\Lambda) \leq \lMdim(\Lambda) \leq \uMdim(\Lambda)
= \delta(G) = \dim(\Lambda_c) \leq \dim(\Lambda) ,$$
so again, all these numbers are equal.  Thus in both cases
the Minkowski dimension exists and equals the Hausdorff dimension.
The second case also proves \thmref{delta=dim}.
$\qed$

Our proof of \thmref{delta=uM} breaks up into three cases:
\begin{enumerate}
\item $\Omega(G)/G$ is compact.
\item $\Omega(G)/G$ has punctures  and $G$ is geometrically finite.
\item $\Omega(G)/G$ has punctures  and $G$ is geometrically 
infinite (i.e., $\dim(\Lambda(G)) =2$).
\end{enumerate}

Note that the third case is not needed for \thmref{H=M}, 
since the two dimensions must agree whenever the Hausdorff
dimension is $2$.

It is worth also noting that our proof uses \thmref{gi->dim=2}, but 
this is  not necessary. It is possible to prove  \thmref{delta=uM}  without
using \thmref{gi->dim=2}. This argument  uses
 an intricate stopping time argument
involving the ``good'' and ``bad'' horoballs and is completely 
``elementary'' (e.g., two dimensional), but using 
\thmref{gi->dim=2} allows us to give a simpler proof. 
It may be  possible that the implication could be reversed and 
\thmref{delta=uM} used to give an alternate proof of 
\thmref{gi->dim=2} but we have not done this (it would suffice
to show that if $G$ is geometrically infinite and $\delta(G)<2$
then $\delta(G) < \uMdim(\Lambda(G))$).

Since $\delta(G) = \dim(\Lambda_c) \leq \dim(\Lambda)$ and 
$\dim(\Lambda) \leq \uMdim(\Lambda)$ we already know that 
$$ \delta(G) \leq \uMdim( \Lambda),$$
but here is a simple proof of this fact which does not 
require \thmref{delta=conical1}.
Chose a point $z_0 \in \Omega(G)$. By the \lemref{lambda=up}
the critical exponent $\delta(G)$  for the Poincar{\'e}
series of $G$ is the same as the critical exponent for the 
series 
$$ \sum_{g\in G} d(g(z))^\alpha.$$
Because each Whitney square contains only a bounded 
number of images of $z_0$, 
$$ \sum_{g\in G} d(g(z_0))^\alpha
\leq C \sum_j \ell(Q_j)^\alpha.$$
Hence the sum on the left converges whenever the sum on the 
right does. Thus 
$$\delta(G) \leq \kappa(\Lambda) \leq \uMdim(\lambda),$$
for any Kleinian group $G$.

We will now start the three cases of \thmref{delta=uM}.

{\it Proof of case (1) in \thmref{delta=uM}:}
Now suppose $G$ is finitely generated, $\Omega(G)/G = R_1\cup
\dots \cup R_N$ is  a union of  compact  surfaces 
and $\Lambda(G)$ has zero area.
Let $\{Q_j\}$ be a Whitney decomposition of $\Omega(G)$. 
By the previous lemma there is an $\epsilon >0$ so that 
every $Q_j$ contains a hyperbolic ball of radius $\epsilon$.

By the compactness of the surfaces $R_1, \dots, R_N$, we can 
choose a finite number of  points $ E=\{ z_1, \dots, 
z_M\} \subset \Omega(G)$, so that $E$ projects to an $\epsilon$-dense
set in $R_1 \cup \dots \cup R_N$, i.e., every point of $\Omega$
is within distance $\epsilon$ of some  point of 
$ G(E) =\cup_{i=1}^M \cup_{g\in G} g(z_i)$.
Thus each Whitney square $Q_j$ contains at least one point of $G(E)$.
Therefore
$$ \sum_j \ell(Q_j)^\alpha \leq C
\sum_{z\in G(E)} d(z)^\alpha.$$
By \lemref{B3=S2}, the infinite series on the 
right hand side  converges for $\alpha > \delta(G)$, hence
so does the left hand side, i.e., $\kappa(\lambda) \leq 
\delta(G)$. Since we already proved the opposite inequality 
we have $\kappa(\Lambda) = \delta(G)$.

If $\Lambda(G)$ has zero area, then  \lemref{kappa=uMdim}
implies 
$$ \uMdim(\Lambda) = \kappa(\Lambda)= \delta(G),$$
as desired.
$\qed$

{\it Proof of case (2) in \thmref{delta=uM}:}
In this case we assume that $G$ is geometrically finite, 
so by \corref{gf->dc} every  parabolic cusp of $\Omega(G)$
is double cusped. Thus by part (2) of \lemref{horoballs}
there is an $\eta>0$ so that every horoball in $\Omega(G)$ 
is $\eta$-good.

Let $\epsilon >0$ and choose $n_0$ so that 
$$ N(\Lambda, 2^{-n_0}) \geq 100\cdot 2^{n_0(D- \epsilon/2)}.$$
Let $r= 2^{-n_0}$. Let $\Cal S = \{ S_k\}$ be a collection of
$2^{n_0(D-\epsilon)}$  squares of size $r$ so that the 
triples $3S_k$ are pairwise disjoint and  $\frac 13 S_k \cap 
\Lambda \ne \emptyset$ for each $k$.

For each  $\eta$-good horoball with $\diam(B) \geq r/3$, let
$\Cal G_B$ be the collection of squares in $\Cal S$ which 
so that $\frac 13 S$ hits $B$. Let $\Cal G$ be the union of 
all the $\Cal G_B$.
Note that the number of squares in $\Cal G_B$ is at 
most $C \diam(B)/r$. Therefore if $R = \diam (B)$, 
\begin{eqnarray*}
 \sum_{S \in \Cal G_B} \ell(S)^{D- \epsilon} 
&\leq&  C \frac Rr r^{D -\epsilon}\\
&=& C R r^{D -1 - \epsilon}\\
&\leq& C R R^{D -1 - \epsilon}\\
&=& C R^{D - \epsilon}
\end{eqnarray*}
 By part (1) of \lemref{horoballs} there
is a orbit point $z \in G(E) \cap B$ such that
$$ d(z) \sim \diam(B).$$
For this point,
 $$d(z)^{D - \epsilon}  \geq C 
 \sum_{S \in \Cal G_B} \ell(S)^{D- \epsilon} .$$

If more than half the squares in $\Cal S$ belong to 
$\Cal G$ then this argument shows 
$$\sum_{z \in G(E)} d(z)^{D- \epsilon} 
\geq \frac 12 C 2^{n_0 \epsilon /2}.$$
If this happens for arbitrarily  large $n_0$ then we have shown 
that $ \delta(G) \geq D - \epsilon$, as desired.

Thus we may assume that fewer than half the elements of $\Cal S$ are in 
$\Cal G$. In this case, part (1) of \lemref{horoballs}
implies that  each horoball $B$ is associated to an orbit point
$z$ so that 
$$ \text{diam}(B) \sim d(z).$$
Moreover, if we sum over all Whitney square hitting $B$ 
then 
$$\sum_{Q\cap B \ne \emptyset} \ell(Q)^\alpha \sim \diam(B)^\alpha
\sim d(z)^\alpha.$$
Therefore the proof in this case is completed exactly as in the 
compact case. $\qed$

{\it Proof of case (3) of \thmref{delta=uM}:}
We can  consider the final case. Suppose that 
$G$ is geometrically infinite (so $\dim(\Lambda(G))=2$
by \thmref{gi->dim=2}) and assume $\text{area}(\Lambda(G)) =0$.
By \lemref{kappa=uMdim} $\kappa (\Lambda(G)) = \uMdim
(\Lambda(G))= \dim(\Lambda(G)) =2$.
Let $\{ B_j\}$ be a listing of the horoballs in $\Omega(G)$.
Each horoball contains a Whitney square of comparable size so 
for any $\epsilon >0$
$$\sum_j \diam(B_j)^{2 - \epsilon/2} = \infty.$$
Thus if 
$$\Cal B_n = \{ B_j: 2^{-n-1} \leq \diam(B_j) < 2^{-n} \},$$
and $N_n = \# \Cal B_n$, we must have 
$$ N_n   \geq 2^{n(2- \epsilon)},$$
for infinitely many values of $n$. Fix a value of $n_0$ 
where this holds and note that for at least half the 
balls  $B$ in $\Cal B_n$  there is a second ball $B'
\in \Cal B_n$ such that
$$\dist(B, B') \leq 2^{n(1-\epsilon)} \leq \diam(B)2^{-n \epsilon},$$
(otherwise we would have so many disjoint  balls that we 
contradict the assumption that $\Lambda(G)$ has diameter $1$).
By part (6) of \lemref{horoballs}, this implies that 
both $B$ and $B'$ are $2^{-3n\epsilon}$-good horoballs. 
Let $\Cal G_{n_0} \subset \Cal B_{n_0}$ be the subcollection 
of $2^{-3n\epsilon}$-good horoballs. For  any $\eta$-good
horoball 
$B$ let  $z$  be the point given in part (1) of 
\lemref{horoballs} and let $h$ be the primitive parabolic 
element  fixing $B_j$. Then an easy calculation shows
$$ \sum_{k \in {\Bbb Z}} d(h^k(z))^\alpha \geq  C \diam(B)^\alpha
\eta^{\alpha -1},$$
for any $0< \alpha < 2$ and some $C$ depending on $G$ and 
$\alpha$.

Thus  if $z_j$ is the good point in $B_j$ given by (1) of 
\lemref{horoballs}, 
\begin{eqnarray*}
 \sum_{ z \in G(E)} d(z)^\alpha 
&\geq& 
\sum_{B_j \in \Cal G_{n_0}} \sum_{k \in {\Bbb Z}} d(h^k(z_j))^\alpha \\
& \geq& C \sum_{B_j \in \Cal G_{n_0}}\diam(B)^\alpha 2^{-3n_0\epsilon(\alpha-1)}\\
& \geq& C 2^{n_0(2- \epsilon)} 2^{-n_0\alpha} 2^{-3n_0\epsilon(\alpha -1)}\\
& \geq& C 2^{n_0[(2- \epsilon)-\alpha-3\epsilon(\alpha -1)]}\\
& \geq& C 2^{n_0[(2 - \alpha)- \epsilon -3 \epsilon(\alpha-1)]}\\
\end{eqnarray*}
If we first choose any $\alpha<2$ and then choose $\epsilon$
sufficiently small, the exponent is positive. Thus 
$\delta(G) > \alpha $ for any $\alpha <2$, i.e., $\delta(G) =2$.
This completes the proof of the third and final case of 
\thmref{delta=uM}.  $\qed$

This finishes the proofs of \thmref{H=M} and \thmref{delta=uM}.

\newpage

\section{Teichm{\"u}ller spaces} \label{Teich}

In this section we shall consider $\dim(\Lambda(G))$ as 
a function on the closure of Teichm{\"u}ller space 
$T(S)$ of a finite type hyperbolic surface $S$.

Given a finite type surface $S$ (compact with a finite
number of punctures (possible none)),
the Teichm{\"u}ller space $T(S)$ is the set of equivalence
classes of quasiconformal mapping of $S$ to itself. 
Each such is represented by a Beltrami differential $\mu$ 
which may be lifted to a Beltrami differential  $\mu$ on the 
the upper halfplane, $\Bbb H$.
Let $\Gamma$ be a Fuchsian group acting on  $\Bbb H$ 
such that ${\Bbb H} /\Gamma = S$.
There is a quasiconformal mapping $F$  of the plane  which 
fixes $0,1,\infty$ and such that $\Dbar F/ \partial F = \mu$
on $\Bbb H$ and so that $F$ is conformal on the lower half plane.
On the lower half-plane the Schwarzian derivative $S(F)$ satisfies
$$ \| S(F)\|= \sup_z |\text{Im}(z)|^2  | S(F)(z)|   \geq  6< \infty.$$
This realizes $T(S)$ as a bounded subset of a Banach space and 
gives a metric on $T(S)$. The closure of $T(S)$ with respect to 
this metric is denoted $\overline{T(S)}$ and the boundary 
by $\partial T(S)$. 
Points of $\overline{T(S)}$ may be identified with 
certain Kleinian groups which are isomorphic to $\Gamma$.
 Moreover, convergence in the Teichm{\"u}ler
metric above implies algebraic convergence of the groups.
Recall that a group is called degenerate if $\Omega(G)$ has 
one component and this component is simply connected. 
Such groups must be geometrically infinite by a result of 
Greenberg \cite{Greenberg66}.
$G \in \partial T(S)$ is called a cusp if there is a hyperbolic 
element in $\Gamma$ which becomes parabolic in $G$.
Bers showed $\partial T(S)$ consists entirely of 
degenerate groups and cusps and that degenerate groups form 
a dense $G_\delta$ set
in $\partial T(S)$ in \cite{Bers70} (in fact, the cusps lie on 
a countable union of real codimension 2 surfaces).
McMullen \cite{McMullen91} proved there is a dense set of 
geometrically finite cusps in $\partial T(S)$.

Recall that \thmref{Hdim=lsc} says that if $\{ g_n\}$ converges 
algebraically to $G$ then  
$$ \liminf_{n\to \infty} \dim(\Lambda(G_n))) \geq \dim(\Lambda (G)).$$
One special case where this holds
is for $ G \in  \overline{T(S)}$, 
the closure of  the Teichm{\"u}ller space of a finite type
hyperbolic Riemann surface $S$. Since $\Lambda(G)$ is at most 
$2$ and is lower semi-continuous, it is continuous whenever it 
takes the value $2$ (i.e., at the geometrically infinite groups).
Since these points are dense on the boundary of Teichm{\"u}ller 
space, this function must be discontinuous at the geometrically finite
cusps on the boundary. Thus 

\begin{cor} \label{dim=lsc,T(S)}
Suppose $S$ is a hyperbolic Riemann surface of finite type.
Then $\dim(\Lambda(G))$ is lower semi-continuous on 
$\overline{T(S)}$,
 and is continuous everywhere except at 
the  geometrically finite cusps in $\partial T(S)$ where it is 
discontinuous.
\end{cor}

This also show that equality in \corref{delta=lsc} and  \thmref{Hdim=lsc}
 need not occur, 
because  a geometrically finite cusp ($\delta(G) < 2$) can be approximated by 
by geometrically finite groups which are close 
to degenerate groups (so $\delta(G_n) \to 2$).
The discontinuity at the geometrically finite cusped had been proved
earlier by Taylor in \cite{Taylor94}; he showed that for each 
geometrically finite cusp $G$ there is a sequence $G_n \to G$ 
algebraically, but $G_n \to H$ geometrically where $H$ is a 
geomtrically finite group containing $G$ and 
$$\dim(\Lambda(G_n)) \to \dim(\Lambda(H)) > \dim(\Lambda(G)).$$

  If $f$ is lower
semi-continuous then $\{ f\leq \alpha\}$ is closed. Thus,

\begin{cor}
The set $ E_\alpha=\{G\in \overline{T(S)}: \dim(\Lambda(G)) 
\leq \alpha \}$ is closed in $\overline{T(S)}$.
The set $ F_\alpha=\{ G \in \partial T(S): \dim(\Lambda(G)) \leq 
\alpha <2\}$ is a closed, nowhere  dense subset of 
$\partial T(S)$.
\end{cor}

Since a lower semi-continuous functions takes a minimum on a compact set,
\thmref{thm1} implies

\begin{cor}\label{dimbound>1}
$\dim(\Lambda(G))$ takes a minimum value on 
$\partial T(S)$ and this minimum is strictly larger than $1$.
\end{cor}

It is not clear where the minimum occurs. Canary has suggested 
it might occur at the cusp group corresponding to shrinking a 
minimum length geodesic on $S$ to a parabolic, since this 
requires the ``least'' deformation of the Fuchian group (in some 
sense). Since $\dim(\Lambda(G))$ takes a minimum on $\partial
T(S)$ which is $>1$, any group in $\overline{T(S)}$ with small 
enough dimension must be quasi-Fuchsian. Thus 
  



\begin{cor}
Suppose $\{G_n\}$ is a sequence of quasiconformal deformations
of a Fuchsian group $G $ (i.e., $\{G_n\}$ is a sequence in 
$T(S)$, $S = {\Bbb D}/ G$). If $\dim(\Lambda(G_n)) \to 1$, then 
$G_n \to G$.
\end{cor}
 
\begin{pf} 
This is immediate from \thmref{dim=lsc} and the fact 
(from \thmref{thm1}) that $G$ is the only point in 
$\overline{T(S)}$ where $\dim (\Lambda)=1$.
\end{pf}

\begin{cor} \label{qcdef}
Suppose $G$ is a finitely generated   and 
has  a simply connected  invariant component $\Omega_0$
(possibly not unique). Let $\Omega_0 /G =S$.
 Then for any $\epsilon >0$ there
is a $\delta$, depending only on  $S$ and $\epsilon$, such that 
$\dim(\Lambda(G)) < 1+\delta$ implies $G$ is a 
$\epsilon$-quasiconformal deformation of a Fuchsian group.
\end{cor}

\begin{pf}
If $G$ is a quasi-Fuchsian group this follows from the 
previous result. Maskit \cite{Maskit70} proved that
a finitely generated Kleinian group with two invariant 
components is quasi-Fuchsian, so we may now assume 
$G$ is a $b$-group (i.e., 
has exactly one simply connected, invariant component).
If $G$ is geometrically infinite then $\dim(\Lambda(G))=2$, 
which contradicts our assumption. Therefore, $G$ must be
geometrically finite. Abikoff \cite{Abikoff75}
 proved that every geometrically 
finite $b$ group covering $S$ is on the boundary of 
the Teichm{\"u}ller space  $T(S)$, and so its dimension is
bounded away from $1$ by \corref{dimbound>1}.
\end{pf}

It is not true that the $\epsilon$ in \corref{qcdef} can be
taken to depend only on the topological type of $S$ (e.g., 
the number of generators of $G$). For example, given a 
surface with punctures $S$, it is possible to use the combination 
theorems to construct a $b$-group $G$ with 
$\dim(\Lambda(G))$ as close to one as we wish and so 
that $\Omega_0/G$ is homeomorphic to (thought not conformally 
equivalent to ) $S$.

Larman showed that  there is an $\epsilon_0$, such that
if $\{D_j\}$ is a collection of three or 
more disjoint open disks then the dimension of 
$\Cbar \setminus \cup_j D_j$ is larger than $1+\epsilon_0$.
A careful  reading of Larman's paper \cite{Larman67} shows that
his proof gives

\begin{prop}
There is an $\epsilon_0 >0$ such that if $\{ D_j\}$ is any 
collection of three or more disjoint  open $\epsilon_0$-quasidisks, then 
$\dim({\Bbb C} \setminus \cup_j D_j ) > 1+ \epsilon_0$.
\end{prop}

Recall that a {\it web group} is a finitely generated 
Kleinian groups each of whose component subgroups is 
quasi-Fuchsian. 
Suppose $G$ is web group. If $\Omega(G)$ has only two 
components then $G$ is quasi-Fuchsian \cite{Maskit88}).
So suppose $G$ has three or more components  and 
let $\{G_1, \dots, G_n\}$ representatives of each conjugacy 
class of component subgroups. By the last corollary either 
one of these has limit set with dimension $> 1+\delta$ or all 
are $\epsilon_0$-quasicircle. In the latter case. Larman's 
theorem implies that $\Lambda(G)$ has dimension bigger than 
$\epsilon_0$. In either case the dimension is bounded away 
from $1$ by a number which only depends on the conformal
structure of $\Omega(G)/G$.
of generators for each component subgroup.


\begin{cor} \label{web}
Suppose $G$ is a finitely generated web group which  is not 
quasi-Fuchsian. Then $\dim(\Lambda ) > 1+ \epsilon $ where 
$\epsilon$  depends  only on the conformal types of the 
 components of  $\Omega(G)/G$.
\end{cor}

It is not true that  the dimension of limit sets of proper web groups
(i.e., not quasi-Fuchsian) is bounded uniformly away from $1$.
Canary, Minsky and Taylor have constructed examples of 
proper web groups (with a fixed number of generators)
whose limits sets have dimension arbitarily close to 1 (personal
communication).

\end{document}